\documentstyle[11pt]{article}
\def\Bbb R{{\rm \bf R}}
\def\proclaim#1{\vskip2mm{\bf #1}\em}
\def\endproclaim{\em \vskip2mm}
\def\tag#1{\eqno(#1)}
\def\gathered{\begin{array}{c}}
\def\endgathered{\end{array}}
\def\text{\mbox}

\begin{document}

\title {On finding a penetrable obstacle using a single electromagnetic wave in the time domain}
\author{Masaru IKEHATA\footnote{
Laboratory of Mathematics,
Graduate School of Advanced Science and Engineering,
Hiroshima University, Higashihiroshima 739-8527, JAPAN}}
\maketitle
 
\begin{abstract}
The time domain enclosure method is one of analytical methods for inverse obstacle problems governed by partial
differential equations in the time domain.
This paper considers the case when the governing equation is given by the Maxwell system
and consists of two parts.  The first part establishes the base of the time domain enclosure method
for the Maxwell system using a single set of the solutions over a finite time interval
for a general (isotropic) inhomogeneous medium in the whole space.  It is a system of asymptotic inequalities
for the indicator function which may enable us to apply the time domain enclosure method to the problem of finding unknown penetrable
obstacles embedded in various background media.  As a first step of its expected applications, the case when 
the background medium is homogeneous and isotropic, is considered and the time domain enclosure method is realized.  
This is the second part.


\noindent
AMS: 35R30, 78A46, 78A40, 35R05, 35Q60, 35B40, 35R25

\end{abstract}


\section{Introduction}

The time domain enclosure method using a single set of observation data is one of analytical methods for inverse obstacle problems governed by 
various partial differential equations in the time domain.  The method tells us how to extract information about
unknown obstacles (so-called discontinuity in non destructive testing) from the time domain data
generated by a single solution of the governing equation over a {\it finite time interval} observed at a place not far a way from the obstacle.
It is considered as a time domain realization of the classical enclosure method \cite{I1, IE}.
The time domain enclosure method using a single set of observation data was initiated by the author in \cite{I4} by considering one space dimensional case for the heat and wave equations
and developed for three space dimensional case in \cite{IE0, IW2, IWa, IWd, IK, IKK} for scalar wave equations in the whole space
or an exterior domain and recently \cite{Iwave, EV} in a bounded domain or \cite{Ithermo, EVI} for the heat equation or
a coupled system in a linear theory of thermoelasticity \cite{C}.

In a series of articles \cite{IMax, IMax2, IMax3} the method has been applied also to inverse obstacle problems
governed by the Maxwell system in an exterior domain.
The unknown obstacles considered therein are impenetrable ones.
Based on these studies the next goal is to consider the problem of finding unknown {\it penetrable obstacles} embedded in various background media
which has various possible applications to non destructive evaluation, radar, microwave tomography and through-the-wall imaging, etc..
The purpose of this paper is to establish the base of the time domain enclosure method using a single set of observation data for the Maxwell system
for penetrable obstacles embedded in the whole space and gives one of possible applications.

Let us describe the mathematical formulation of the problem.

Let $0<T<\infty$.
We denote by $\mbox{\boldmath $E$}=\mbox{\boldmath $E$}(x,t)$ and $\mbox{\boldmath $H$}=\mbox{\boldmath $H$}(x,t)$, $(x,t)\in\Bbb R^3\times\,[0,\,T]$ the electric field and the magnetic field, respectively;
the functions $\epsilon=\epsilon(x)$, $\mu=\mu(x)$ and $\sigma=\sigma(x)$
denote the electric permittivity, 
the magnetic permeability 
and electrical conductivity of an isotropic medium occupying the whole space assumed to be essentially bounded functions on $\Bbb R^3$
such that $\text{ess}.\,\inf_{x\in\Bbb R^3}\epsilon(x)>0$, $\text{ess}.\,\inf_{x\in\Bbb R^3}\mu(x)>0$
and $\text{ess}.\,\inf_{x\in\Bbb R^3}\sigma(x)\ge 0$, respectively.

Let $D$ be a nonempty bounded open subset of $\Bbb R^3$ such that
$\Bbb R^3\setminus\overline D$ is connected. 
We assume that $\epsilon$, $\mu$ and $\sigma$ take the form
$$\displaystyle
\epsilon_r(x)\equiv\frac{\epsilon(x)}{\epsilon_0(x)}=\left\{
\begin{array}{ll}
\displaystyle
1, &
\displaystyle
x\in \Bbb R^3\setminus D,\\
\\
\displaystyle
1+e(x), &
\displaystyle
x\in D,
\end{array}
\right.
$$
$$\displaystyle
\mu_r(x)\equiv \frac{\mu(x)}{\mu_0(x)}
=\left\{
\begin{array}{ll}
\displaystyle
1, &
\displaystyle
x\in \Bbb R^3\setminus D,\\
\\
\displaystyle
1+m(x), & 
\displaystyle
x\in D,
\end{array}
\right.
$$
and
$$\displaystyle
\sigma(x)
=\left\{
\begin{array}{ll}
\displaystyle
\sigma_0(x), & 
\displaystyle
x\in \Bbb R^3\setminus D,\\
\\
\displaystyle
\sigma_0(x)+h(x), & 
\displaystyle
x\in D,
\end{array}
\right.
$$
where $\epsilon_0=\epsilon_0(x)$, $\mu_0=\mu_0(x)$ and $\sigma_0=\sigma_0(x)$ denote the electric permittivity,
magnetic permeability and electric conductivity of the background medium occupying $\Bbb R^3$
and assumed to be essentially bounded functions on $\Bbb R^3$
such that $\text{ess}.\,\inf_{x\in\Bbb R^3}\epsilon_0(x)>0$, $\text{ess}.\,\inf_{x\in\Bbb R^3}\mu_0(x)>0$
and $\text{ess}.\,\inf_{x\in\Bbb R^3}\sigma_0(x)\ge 0$;
$e=e(x)$, $m=m(x)$ and $h=h(x)$ are essentially bounded functions on $D$.

The set $D$ is a mathematical model of a set of penetrable obstacles embedded in the background medium.
Note that at this stage, we never assume any regularity of $\partial D$, $\mu_0$, $\epsilon_0$,
$\sigma_0$, $\mu$, $\epsilon$ and $\sigma$ nor jump conditions on $\mu_r$, $\epsilon_r$
across $\partial D$.

We assume that the fields $\mbox{\boldmath $E$}$ and $\mbox{\boldmath $H$}$
are induced only by the current density $\mbox{\boldmath $J$}=\mbox{\boldmath $J$}(x,t)$ at $t=0$  whose support is localized 
outside $D$.
The governing equations of fields $\mbox{\boldmath $E$}$ and $\mbox{\boldmath $H$}$
take the form
$$
\left\{
\begin{array}{ll}
\displaystyle
\epsilon\frac{\partial\mbox{\boldmath $E$}}{\partial t}
-\nabla\times\mbox{\boldmath $H$}=-\sigma\mbox{\boldmath $E$}+\mbox{\boldmath $J$},
& 
\displaystyle
(x,t)\in \Bbb R^3\times\,]0,\,T[,
\\
\\
\displaystyle
\mu\frac{\partial\mbox{\boldmath $H$}}{\partial t}
+\nabla\times\mbox{\boldmath $E$}=\mbox{\boldmath $0$},
&
\displaystyle
(x,t)\in \Bbb R^3\times\,]0,\,T[,
\\
\\
\displaystyle
\mbox{\boldmath $E$}\vert_{t=0}=\mbox{\boldmath $0$},
&
\displaystyle
x\in\Bbb R^3,
\\
\\
\displaystyle
\mbox{\boldmath $H$}\vert_{t=0}=\mbox{\boldmath $0$},
&
\displaystyle
x\in\Bbb R^3.
\end{array}
\right.
\tag {1.1}
$$

Fix a large (to be determined later) $T<\infty$.
Let $B$ be the open ball centered at a point $p$ with {\it very small} radius $\eta$
and satisfy $\overline B\cap\overline D=\emptyset$.
There are several choices of the current density $\mbox{\boldmath $J$}$
as a model of the antenna.
In this paper, as considered before in \cite{IMax, IMax2, IMax3} we assume that $\mbox{\boldmath $J$}$
takes the form
$$\displaystyle
\mbox{\boldmath $J$}(x,t)
=f(t)\chi_B(x)\mbox{\boldmath $a$},
$$
where $\mbox{\boldmath $a$}\not=\mbox{\boldmath $0$}$ is a constant unit vector,
$\chi_B$ denotes the characteristic function of $B$ and $f\in
H^1(0,\,T)$ with $f(0)=0$. Note that $\chi_B(x)$ has discontinuity
across the sphere $\partial B$.

Now we are ready to state our problem.

{\bf Problem.}  Fix a large $T$ (to be specified later).  Generate fields $\mbox{\boldmath $E$}$ and $\mbox{\boldmath $H$}$ by $\mbox{\boldmath $J$}$ and
observe $\mbox{\boldmath $E$}$ on $B$ over time interval $]0,\,T[$.
Extract information about the geometry of $D$ from the observed data.

To describe the results for the problem we introduce  an indicator function in the {\it time domain enclosure method}.
For the purpose and the use in other parts we introduce some fields.

Let $\tau>0$ and set
$$
\left\{
\begin{array}{ll}
\displaystyle
\mbox{\boldmath $W$}_e=\mbox{\boldmath $W$}_e(x,\tau)=\int_0^Te^{-\tau t}\mbox{\boldmath $E$}(x,t)dt,
&
\displaystyle
x\in\Bbb R^3,
\\
\\
\displaystyle
\mbox{\boldmath $W$}_m=\mbox{\boldmath $W$}_m(x,\tau)=\int_0^Te^{-\tau t}\mbox{\boldmath $H$}(x,t)dt,
&
\displaystyle
x\in\Bbb R^3.
\end{array}
\right.
\tag {1.2}
$$
Let $\mbox{\boldmath $E$}_0=\mbox{\boldmath $E$}_0(x,t)$ and $\mbox{\boldmath $H$}_0=\mbox{\boldmath $H$}_0(x,t)$ satisfy
$$\displaystyle
\left\{
\begin{array}{ll}
\displaystyle
\epsilon_0\frac{\partial\mbox{\boldmath $E$}_0}{\partial t}
-\nabla\times\mbox{\boldmath $H$}_0=-\sigma_0\mbox{\boldmath $E$}_0+\mbox{\boldmath $J$},
&
\displaystyle
(x,t)\in\Bbb R^3\times\,]0,\,T[,
\\
\\
\displaystyle
\mu_0\frac{\partial\mbox{\boldmath $H$}_0}{\partial t}
+\nabla\times\mbox{\boldmath $E$}_0=\mbox{\boldmath $0$},
&
\displaystyle
(x,t)\in\Bbb R^3\times\,]0,\,T[,
\\
\\
\displaystyle
\mbox{\boldmath $E$}_0\vert_{t=0}=\mbox{\boldmath $0$},
&
\displaystyle
x\in\Bbb R^3,
\\
\\
\displaystyle
\mbox{\boldmath $H$}_0\vert_{t=0}=\mbox{\boldmath $0$},
&
\displaystyle
x\in\Bbb R^3
\end{array}
\right.
\tag {1.3}
$$
and set
$$
\left\{
\begin{array}{ll}
\displaystyle
\mbox{\boldmath $V$}_e=\mbox{\boldmath $V$}_e(x,\tau)=\int_0^Te^{-\tau t}\mbox{\boldmath $E$}_0(x,t)dt,
&
\displaystyle
x\in\Bbb R^3,
\\
\\
\displaystyle
\mbox{\boldmath $V$}_m=\mbox{\boldmath $V$}_m(x,\tau)=\int_0^Te^{-\tau t}\mbox{\boldmath $H$}_0(x,t)dt,
&
\displaystyle
x\in\Bbb R^3.
\end{array}
\right.
\tag {1.4}
$$

{\bf\noindent Definition 1.1.}  Define the indicator function of the time domain enclosure method
by the formula
$$\begin{array}{ll}
\displaystyle
I(\tau)=\int_{\Bbb R^3}\mbox{\boldmath $f$}\cdot(\mbox{\boldmath $W$}_e-\mbox{\boldmath $V$}_e)dx,
&
\displaystyle
\tau>0,
\end{array}
$$
where
$$\displaystyle
\mbox{\boldmath $f$}=\mbox{\boldmath $f$}(x,\tau)
=\int_0^Te^{-\tau t}\mbox{\boldmath $J$}(x,t)dt
=\tilde{f}(\tau)\chi_B(x)\mbox{\boldmath $a$}
$$
and
$$\displaystyle
\tilde{f}(\tau)=\int_0^T e^{-\tau t}f(t)dt.
$$
It is easy to see that we have, as $\tau\longrightarrow\infty$
$$\displaystyle
\Vert\mbox{\boldmath $f$}(\,\cdot\,,\tau)\Vert_{L^2(\Bbb R^3)}=O(\vert\tilde{f}(\tau)\vert)=O(\tau^{-3/2}).
\tag {1.5}
$$
Note that the support of $\mbox{\boldmath $f$}(\,\cdot\,,\tau)$ is contained in $\overline B$ and thus the indicator function
can be calculated by using the data $\mbox{\boldmath $E$}(x,t)\cdot\mbox{\boldmath $a$}$, $(x,t)\in B\times\,]0,\,T[$.
Note that this type of indicator function was firstly appeared in \cite{IW2} for scalar wave equations. 

The first result is concerned with upper and lower asymptotic estimates of the indicator function 
in terms of the energy integrals of $\mbox{\boldmath $V$}_e$ and $\mbox{\boldmath $V$}_m$ given by (1.4) over $D$.

\proclaim{\noindent Theorem 1.1.}  Fix an arbitrary $T$.
We have, as $\tau\longrightarrow\infty$
$$
\displaystyle
I(\tau)
\le 
\tau\int_{\Bbb R^3}
\left\{
\frac{\tilde{\epsilon_0}}{\tilde{\epsilon}}(\tilde{\epsilon_0}-\tilde{\epsilon})
\vert\mbox{\boldmath $V$}_e\vert^2dx
+
(\mu-\mu_0)
\vert\mbox{\boldmath $V$}_m\vert^2
\right\}dx
+O(\tau^{-5/2}e^{-\tau T})
\tag {1.6}
$$
and
$$\displaystyle
I(\tau)
\ge 
\tau\int_{\Bbb R^3}\left\{(\tilde{\epsilon_0}-\tilde{\epsilon})\vert\mbox{\boldmath $V$}_e\vert^2
+\frac{\mu_0}{\mu}(\mu-\mu_0)\vert\mbox{\boldmath $V$}_m\vert^2\right\}dx
+O(\tau^{-5/2}e^{-\tau T}),
\tag {1.7}
$$
where
$$\displaystyle
\tilde{\epsilon}=\epsilon+\frac{\sigma}{\tau}
\tag {1.8}
$$
and
$$\displaystyle
\tilde{\epsilon_0}=\epsilon_0+\frac{\sigma_0}{\tau}.
\tag {1.9}
$$

\endproclaim

Note that we have
$$\left\{
\begin{array}{l}
\displaystyle
(\tilde{\epsilon_0}-\tilde{\epsilon})+\frac{(\tilde{\epsilon_0}-\tilde{\epsilon})^2}{\tilde{\epsilon}}=
\frac{\tilde{\epsilon}_0}{\tilde{\epsilon}}\,(\tilde{\epsilon_0}-\tilde{\epsilon}),
\\
\\
\displaystyle
\frac{\mu_0}{\mu}(\mu-\mu_0)+\frac{(\mu-\mu_0)^2}{\mu}=\mu-\mu_0.
\end{array}
\right.
\tag {1.10}
$$
Thus the first term of the right-hand side on (1.6) is greater than
that of the right-hand side on (1.7).

One can expect that Theorem 1.1 enables us to apply the time domain enclosure method to give a solution
to Problem.  
Here, as a first step we consider the case when the background medium is {\it homogeneous},
that is, $\epsilon_0$, $\mu_0$ and $\sigma_0$ are {\it constant}.

We introduce the following jump conditions (A.I) and (A.II) for $\mu_r$ and $\epsilon_r$ on the unknown
obstacle $D$:

(A.I)  there exists a positive constant $C_1$ such that, for almost all $x\in D$
$$\displaystyle
\left(1-\frac{1}{\epsilon_r(x)}\right)+(1-\mu_r(x))\ge C_1.
\tag {1.11}
$$

(A.II) there exists a positive constant $C_2$ such that, for almost all $x\in D$
$$\displaystyle
(1-\epsilon_r(x))+\left(1-\frac{1}{\mu_r(x)}\right)\ge C_2.
\tag {1.12}
$$

{\bf\noindent Remark 1.1.}
If $\mu_r$ satisfies, for almost all $x\in D$
$$\displaystyle
\mu_r(x)\le 1,
$$
and there exists a positive constant $C_3$ such that, for almost all $x\in D$
$$\displaystyle
\epsilon_r(x)\ge 1+C_3,
$$
then (A.I) is satisfied.  Note that this is one of examples of the pair $(\mu_r,\epsilon_r)$ satisfying (A.I) since the restriction (1.11) is a requirement
for the pair not for $\mu_r$ and $\epsilon_r$ independently.  The same remark works also for (1.12) and (A.II).
For this define two regions $A_1$ and $A_2$ in the plane:
$$\begin{array}{l}
\displaystyle
A_1=\left\{(x,y)\,\vert\,y>\frac{1}{2-x},\, 0<x<2\,\right\},
\\
\\
\displaystyle
A_2=\left\{(x,y)\,\vert\, 0<y<2-\frac{1}{x},\,x>\frac{1}{2}\,\right\}.
\end{array}
$$
The sets $A_1$ and $A_2$ are symmetric with respect to the line $y=x$.
We have $\Phi(A_1)=A_2$ and $\Phi(A_2)=A_1$, where
$$\displaystyle
\Phi:(x,y)\longmapsto (\frac{1}{x},\frac{1}{y}).
$$
Roughly speaking, if the pair $(\mu_r,\epsilon_r)$ satisfies condition (A.I), it means that $(\mu_r,\epsilon_r)\in A_1$ and (A.II) means that $(\mu_r,\epsilon_r)\in A_2$.
Note that both $A_1$ and $A_2$ are unbounded, $A_1\cap A_2=\emptyset$ and $\partial A_1\cap\partial A_2=\{(1,1)\}$.

Next we introduce the following conditions (B.I), (B.II) and (B.III) which are concerned with the shape of the surface 
of unknown $D$ or the direction of $\mbox{\boldmath $a$}$ in $\mbox{\boldmath $J$}$ relative to $\partial D$.

Before doing so, we describe some notion from the differential geometry \cite{O} and related elementary facts, see also \cite{IMax}.
Let $\partial D$ be $C^2$ and for each $q\in\partial D$ we denote by $\mbox{\boldmath $\nu$}_q$ the unit outer normal to $\partial D$ at $q$.
Let $\mbox{\boldmath $S$}_q(\partial D)$ denote the {\it shape operator} of $\partial D$ at $q\in\partial D$
with respect to $\mbox{\boldmath $\nu$}_q$.
Recall that the point $p$ is the center of the ball $B$ with $\overline B\cap\overline D=\emptyset$.

The first reflector from the point $p$ is the set of all $q\in\partial D$ such that
$\vert q-p\vert=\text{dist}(\{p\},\partial D)$ and has the expression
$\partial D\cap S(p,\partial D)$, where $S(p,\partial D)=\{x\in\Bbb R^3\,\vert\,\vert x-p\vert=\text{dist}(\{p\},\partial D)\,\}$.
Then for each point $q$ in the first reflector from point $p$ the vector $\mbox{\boldmath $\nu$}_q$ coincides with
the unit inner normal to the sphere $S(p,\partial D)$ at $q$ and the tangent planes of both $\partial D$ and $S(p,\partial D)$ at
$q$ are the same.
We denote by $\mbox{\boldmath $S$}_q(S(p,\partial D))$ the shape operator of sphere $S(p,\partial D)$ at $q$ with respect to $\mbox{\boldmath $\nu$}_q$
which is the outer normal to $\partial D$ at $q$.
It is known that $\mbox{\boldmath $S$}_q(S(p,\partial D))-\mbox{\boldmath $S$}_q(\partial D)\ge 0$ as the quadratic form on the common tangent space at
each point $q$ in the first reflector from point $p$ and we have the expression
$$\displaystyle
\text{det}\,(\mbox{\boldmath $S$}_q(S(p,\partial D))-\mbox{\boldmath $S$}_q(\partial D))
=\lambda^2-2\lambda H_{\partial D}(q)+K_{\partial D}(q),
$$
where $H_{\partial D}(q)$ is the mean curvature of $\partial D$ at $q$ with respect to $\mbox{\boldmath $\nu$}_q$,
$K_{\partial D}(q)$ the Gauss curvature and $\lambda=(\text{dist}(\{p\},\partial D))^{-1}$.

Now we are ready to state three conditions (B.I), (B.II) and (B.III).

(B.I)  the first reflector from the point $p$ is finite and that,
each point $q$ in the first reflector 
$$\displaystyle
K_{\partial D}(q)>0.
$$

(B.II) the first reflector from the point $p$ is finite and that,
at each point $q$ in the first reflector
$$\displaystyle
\text{det}\,(\mbox{\boldmath $S$}_q(S(p,\partial D))-\mbox{\boldmath $S$}_q(\partial D))>0.
$$

(B.III) there exists a point $q$ in the first reflector from the point $p$ such that 
$\mbox{\boldmath $a$}\times\mbox{\boldmath $\nu_q$}
\not=\mbox{\boldmath $0$}$.

Now we are ready to state a solution to Problem for a penetrable obstacle embedded in a homogeneous background medium.

\proclaim{\noindent Theorem 1.2.}

(i)  If $T<2\sqrt{\mu_0\epsilon_0}\,\text{dist}\,(D,B)$, then we have
$$\displaystyle
\displaystyle
\lim_{\tau\rightarrow\infty}e^{\tau T}I(\tau)=0.
$$

(ii)  Let $f(t)$ satisfy
$$\begin{array}{ll}
\displaystyle
\exists\gamma\in\Bbb R
&
\displaystyle
\liminf_{\tau\rightarrow\infty}\tau^{\gamma}\vert \tilde{f}(\tau)\vert>0.
\end{array}
\tag {1.13}
$$ 
Assume that $\partial D$ is $C^2$ and that the one of the conditions (B.I), (B.II) and (B.II) is satisfied:

Then, we have
$$
\displaystyle
\lim_{\tau\rightarrow\infty}e^{\tau T}I(\tau)
=
\left\{
\begin{array}{ll}
\displaystyle
-\infty,
&
\text{if $T>2\sqrt{\mu_0\epsilon_0}\,\text{dist}\,(D,B)$ and (A.I) is satisfied,}\\
\\
\displaystyle
\infty,
&
\text{if $T>2\sqrt{\mu_0\epsilon_0}\,\text{dist}\,(D,B)$ and (A.II) is satisfied.}
\end{array}
\right.
$$
Besides, if one of (A.I) or (A.II) is satisfied, then we have,
for all $T>2\sqrt{\mu_0\epsilon_0}\,\text{dist}\,(D,B)$
$$
\displaystyle
\lim_{\tau\rightarrow\infty}\frac{1}{\tau}
\log\left\vert
I(\tau)\right\vert
=-2\sqrt{\mu_0\epsilon_0}\,\text{dist}\,(D,B).
\tag {1.14}
$$

\endproclaim

Some remarks are in order.

$\bullet$  There is no assumption on $\sigma$ except for the assumption $\sigma(x)=\sigma_0$ for almost all $x\in\Bbb R^3\setminus D$.

$\bullet$  It is not assumed that $\mu_r\equiv 1$, that is, $\mu(x)=\mu_0$ for almost all $x\in D$.  This is an advantage
of the time domain enclosure method.   The factorization method in the frequency domain in \cite{KG}
which employs infinitely many output corresponding to infinitely many input as the observation data, does not cover 
the case when $\mu_r\not\equiv 1$.  And note that the author does not know whether a time domain factorization method for the Maxwell system
that covers such case exists or not.

$\bullet$  Roughly speaking, the condition (1.13) means that the source at $t=0$ never vanish at infinite order.

$\bullet$  By the {\it signature} of the quantity $e^{\tau T}I(\tau)$ for a sufficiently large $\tau$, one can distinguish,
roughly speaking, whether
$(\mu_r,\epsilon_r)\in A_1$ or $(\mu_r,\epsilon_r)\in A_2$ provided $(\mu_r,\epsilon_r)\in A_1\cup A_2$
and $T$ is sufficiently large.  This is an extraction of a qualitative information about obstacle $D$.

$\bullet$ Clearly the condition (B.I) in the statement (ii) ensures that, at each $q\in\partial D\cap S(p,\partial D)$ there exists an open ball $B'$
centered at $q$ and a positive number $\lambda>0$ such that the set $D\cap B'$ is {\it contained} in 
the open ball $B''=\{x\in\Bbb R^3\,\vert\,\vert x-(q-\lambda\mbox{\boldmath $\nu$}_q)\vert<\lambda\,\}$.  Note that $-\mbox{\boldmath $\nu$}_q$
is the unit inward normal to $\partial B''$ at $q$.  We use only this fact in the proof of (ii) for the case that (B.I) is satisfied.

$\bullet$  The condition (B.II) is a generalization of (B.I).
Note that as $\text{dist}\,(\{p\},\partial D)\rightarrow\infty$ 
we have
$$\displaystyle
\text{det}\,(\mbox{\boldmath $S$}_q(S(p,\partial D))-\mbox{\boldmath $S$}_q(\partial D))
\rightarrow K_{\partial D}(q).
$$

$\bullet$  The condition (B.III) that is concerned with the direction $\mbox{\boldmath $a$}$ relative to $\partial D$ is clearly generic.
Besides, if there exist three points in the first reflector from point $p$, the condition is satisfied.

It should be emphasized that conditions (B.I) and (B.II) can not be controlled by an observer since the obstacle is unknown.
However, condition (B.III) can be resolved by using arbitrary linearly independent directions $\mbox{\boldmath $a$}_j$, $j=1,2$
as done in \cite{IMax2}.  The point is the simple fact: for an arbitrary point $q$ in the first reflector from point $p$
we have $\mbox{\boldmath $a$}_1\times\mbox{\boldmath $\nu$}_q\not=\mbox{\boldmath $0$}$ or
$\mbox{\boldmath $a$}_2\times\mbox{\boldmath $\nu$}_q\not=\mbox{\boldmath $0$}$.

For the purpose define another indicator function by the formula
$$\displaystyle
\mbox{\boldmath $I$}(\tau)
=\sum_{j=1}^2I_j(\tau),
$$
where $I_j(\tau)=I(\tau)$ with $\mbox{\boldmath $J$}$ for $\mbox{\boldmath $a$}=\mbox{\boldmath $a$}_j$, $j=1,2$.

As an easy consequence of a slight modification of the proof of the statement (ii) under condition (B.III) in Theorem 1.2  
we obtain the following corollary.

\proclaim{\noindent Corollary 1.1.}
Assume that $\partial D$ is $C^2$.
Let $f(t)$ satisfy (1.13).
Then, we have:
$$
\displaystyle
\lim_{\tau\rightarrow\infty}e^{\tau T}\mbox{\boldmath $I$}(\tau)
=
\left\{
\begin{array}{ll}
\displaystyle
0,
& \text{if $T<2\sqrt{\mu_0\epsilon_0}\,\text{dist}\,(D,B)$,}\\
\\
\displaystyle
-\infty,
&
\text{if $T>2\sqrt{\mu_0\epsilon_0}\,\text{dist}\,(D,B)$ and (A.I) is satisfied,}\\
\\
\displaystyle
\infty,
&
\text{if $T>2\sqrt{\mu_0\epsilon_0}\,\text{dist}\,(D,B)$ and (A.II) is satisfied.}
\end{array}
\right.
$$

Besides, if one of (A.I) or (A.II) is satisfied, then we have,
for all $T>2\sqrt{\mu_0\epsilon_0}\,\text{dist}\,(D,B)$
$$
\displaystyle
\lim_{\tau\rightarrow\infty}\frac{1}{\tau}
\log\left\vert
\mbox{\boldmath $I$}(\tau)\right\vert
=-2\sqrt{\mu_0\epsilon_0}\,\text{dist}\,(D,B).
\tag {1.15}
$$

\endproclaim

This is an extension of Theorem 1.1 in \cite{IMax2}.
Since we have $\text{dist}(D,B)=\text{dist}(\{p\},\partial D)-\eta$, from the formula (1.15), one gets the quantity
$\text{dist}(\{p\},\partial D)$ and thus the sphere $S(p,\partial D)$.
After knowing the sphere for a fixed $p$,
moving {\it infinitely many} small $B'$ instead of $B$ around $p$ and using (1.15) for $B=B'$, one can determine all the points 
$q\in\partial D\cap S(p,\partial D)$.  This is a typical application of (1.15) type formula in the enclosure method.
See Corollary 1 in \cite{IMax} for the detail.

Finally, we note that the class of solutions of the systems (1.1) and (1.3) 
is just same as those in the articles \cite{NV, NVo} which treat a more complicated model and is based on the framework of \cite{DL}.
We will refrain from giving the details, however the calculation here is also formally understandable.
See also \cite{IW2} in which an exact calculation for scalar wave equations based on the framework of \cite{DL} has been given.

This paper is organized as follows. Theorems 1.1
is proved in Section 2. The proof starts with establishing a representation formula of the indicator function.
This is Proposition 2.1.
Having some energy estimates
which are proved in Lemmas 2.1 and 2.2, we finish the proof of Theorem 1.1.

The proof of Theorem 1.2 is given in Section 3.  The key is to replace $\mbox{\boldmath $V$}_e$ and $\mbox{\boldmath $V$}_m$ in Theorem 1.1
with $\mbox{\boldmath $V$}_e^0$ and $\mbox{\boldmath $V$}_m^0$, respectively which are solutions of the system (3.1)
with large parameter $\tau$.  The system can be derived by formally setting $\epsilon_0\mbox{\boldmath $E$}_0(x,T)=
\mu_0\mbox{\boldmath $H$}_0(x,T)=\mbox{\boldmath $0$}$
in the governing system (2.2) of $\mbox{\boldmath $V$}_e$ and $\mbox{\boldmath $V$}_m$.
The solutions of system (3.1) have been already explicitly constructed in \cite{IMax} in the case when 
$\mu_0$ and $\epsilon_0$ are constant and $\sigma_0=0$ and the construction works also for a constant $\sigma_0\not=0$.
This yields Lemma 3.1 which gives us the explicit representation formulae
of $\vert\mbox{\boldmath $V$}_e^0\vert^2$ and $\vert\mbox{\boldmath $V$}_m^0\vert^2$ outside the support of the source.
By carefully studying their asymptotic behaviour, we obtain the desired results
with the help of Lemmas 3.2 and 3.3 which clarify the meaning of (A.I) and (A.II) in the proof of Theorem 1.2.

In Section 4 some of problems to be solved in the framework of the time domain enclosure method are described.

\section{Proof of Theorem 1.1}

From (1.1) and (1.2) we have
$$
\left\{
\begin{array}{l}
\displaystyle
\nabla\times\mbox{\boldmath $W$}_m-\tau\tilde{\epsilon}\mbox{\boldmath $W$}_e+\mbox{\boldmath $f$}
=e^{-\tau T}\epsilon\mbox{\boldmath $E$}(x,T),
\\
\\
\displaystyle
\nabla\times\mbox{\boldmath $W$}_e+\tau\mu\mbox{\boldmath $W$}_m
=-e^{-\tau T}\mu\mbox{\boldmath $H$}(x,T).
\end{array}
\right.
\tag {2.1}
$$
And also from (1.3) and (1.4) we see that $\mbox{\boldmath $V$}_e$ and $\mbox{\boldmath $V$}_m$ satisfy
$$
\left\{
\begin{array}{l}
\displaystyle
\nabla\times\mbox{\boldmath $V$}_m-\tau\tilde{\epsilon_0}\mbox{\boldmath $V$}_e+\mbox{\boldmath $f$}
=e^{-\tau T}\epsilon_0\mbox{\boldmath $E$}_0(x,T),
\\
\\
\displaystyle
\nabla\times\mbox{\boldmath $V$}_e+\tau\mu_0\mbox{\boldmath $V$}_m
=-e^{-\tau T}\mu_0\mbox{\boldmath $H$}_0(x,T).
\end{array}
\right.
\tag {2.2}
$$
Define
$$
\left\{
\begin{array}{l}
\displaystyle
\mbox{\boldmath $R$}_e=\mbox{\boldmath $W$}_e-\mbox{\boldmath $V$}_e,
\\
\\
\displaystyle
\mbox{\boldmath $R$}_m=\mbox{\boldmath $W$}_m-\mbox{\boldmath $V$}_m.
\end{array}
\right.
\tag {2.3}
$$
From (2.1), (2.2) and (2.3) we see that $\mbox{\boldmath $R$}_e$ and $\mbox{\boldmath $R$}_m$ satisfy
$$
\displaystyle
\nabla\times\mbox{\boldmath $R$}_m-\tau\tilde{\epsilon}\mbox{\boldmath $R$}_e+\tau(\tilde{\epsilon_0}-\tilde{\epsilon})\mbox{\boldmath $V$}_e
=e^{-\tau T}\mbox{\boldmath $F$}
\tag {2.4}
$$
and
$$
\displaystyle
\nabla\times\mbox{\boldmath $R$}_e+\tau\mu\mbox{\boldmath $R$}_m-\tau(\mu_0-\mu)\mbox{\boldmath $V$}_m
=-e^{-\tau T}\mbox{\boldmath $G$},
\tag {2.5}
$$
where
$$
\left\{
\begin{array}{l}
\displaystyle
\mbox{\boldmath $F$}=\epsilon\mbox{\boldmath $E$}(x,T)-\epsilon_0\mbox{\boldmath $E$}_0(x,T),
\\
\\
\displaystyle
\mbox{\boldmath $G$}=\mu\mbox{\boldmath $H$}(x,T)-\mu_0\mbox{\boldmath $H$}_0(x,T).
\end{array}
\right.
\tag {2.6}
$$

First we describe a representation formula of the indicator function.
\proclaim{\noindent Proposition 2.1.}
We have
$$\begin{array}{l}
\displaystyle
\,\,\,\,\,\,
\int_{\Bbb R^3}\mbox{\boldmath $f$}\cdot(\mbox{\boldmath $W$}_e-\mbox{\boldmath $V$}_e)dx\\
\\
\displaystyle
=\tau\int_{\Bbb R^3}\left\{(\tilde{\epsilon_0}-\tilde{\epsilon})\vert\mbox{\boldmath $V$}_e\vert^2
+\frac{\mu_0}{\mu}(\mu-\mu_0)\vert\mbox{\boldmath $V$}_m\vert^2\right\}dx\\
\\
\displaystyle
\,\,\,
+\tau\int_{\Bbb R^3}\left\{\tilde{\epsilon}\vert\mbox{\boldmath $R$}_e\vert^2
+\mu\left\vert\mbox{\boldmath $R$}_m+\frac{\mu-\mu_0}{\mu}\mbox{\boldmath $V$}_m\right\vert^2\right\}dx
\\
\\
\displaystyle
\,\,\,
+e^{-\tau T}
\int_{\Bbb R^3}
\left(\epsilon_0\mbox{\boldmath $E$}_0(x,T)\cdot\mbox{\boldmath $W$}_e
-\epsilon\mbox{\boldmath $E$}(x,T)\cdot\mbox{\boldmath $V$}_e
+\mu \mbox{\boldmath $H$}(x,T)\cdot\mbox{\boldmath $V$}_m
-\mu_0 \mbox{\boldmath $H$}_0(x,T)\cdot\mbox{\boldmath $W$}_m\right)dx
\\
\\
\displaystyle
\,\,\,
+e^{-\tau T}\int_{\Bbb R^3}(\mbox{\boldmath $F$}\cdot\mbox{\boldmath $R$}_e
+\mbox{\boldmath $G$}\cdot\mbox{\boldmath $R$}_m)dx.
\end{array}
\tag {2.7}
$$

\endproclaim

{\it\noindent Proof.}
Taking the inner product of the both sides on (2.4) with $\mbox{\boldmath $R$}_e$, we obtain
$$\displaystyle
(\nabla\times\mbox{\boldmath $R$}_m)\cdot\mbox{\boldmath $R$}_e-\tau\tilde{\epsilon}\vert\mbox{\boldmath $R$}_e\vert^2
+\tau(\tilde{\epsilon_0}-\tilde{\epsilon})\mbox{\boldmath $V$}_e\cdot\mbox{\boldmath $R$}_e
=e^{-\tau T}\mbox{\boldmath $F$}\cdot\mbox{\boldmath $R$}_e,
$$
that is,
$$\displaystyle
\tau(\tilde{\epsilon_0}-\tilde{\epsilon})\mbox{\boldmath $V$}_e\cdot\mbox{\boldmath $R$}_e
=\tau\tilde{\epsilon}\vert\mbox{\boldmath $R$}_e\vert^2-(\nabla\times\mbox{\boldmath $R$}_m)\cdot\mbox{\boldmath $R$}_e
+e^{-\tau T}\mbox{\boldmath $F$}\cdot\mbox{\boldmath $R$}_e.
\tag {2.8}
$$
Next, taking the inner product of the both sides on (2.5) with $\mbox{\boldmath $R$}_m$, we obtain
$$\displaystyle
(\nabla\times\mbox{\boldmath $R$}_e)\cdot\mbox{\boldmath $R$}_m+\tau\mu\vert\mbox{\boldmath $R$}_m\vert^2
-\tau(\mu_0-\mu)\mbox{\boldmath $V$}_m\cdot\mbox{\boldmath $R$}_m
=-e^{-\tau T}\mbox{\boldmath $G$}\cdot\mbox{\boldmath $R$}_m,
$$
that is
$$\displaystyle
\tau(\mu_0-\mu)\mbox{\boldmath $V$}_m\cdot\mbox{\boldmath $R$}_m
=\tau\mu\vert\mbox{\boldmath $R$}_m\vert^2+(\nabla\times\mbox{\boldmath $R$}_e)\cdot\mbox{\boldmath $R$}_m
+e^{-\tau T}\mbox{\boldmath $G$}\cdot\mbox{\boldmath $R$}_m.
\tag {2.9}
$$
Since we have
$$\displaystyle
\nabla\cdot(\mbox{\boldmath $R$}_e\times\mbox{\boldmath $R$}_m)
=\mbox{\boldmath $R$}_m\cdot(\nabla\times\mbox{\boldmath $R$}_e)
-\mbox{\boldmath $R$}_e\cdot(\nabla\times\mbox{\boldmath $R$}_m),
$$
adding (2.8) and (2.9), we obtain
$$\begin{array}{l}
\displaystyle
\,\,\,\,\,\,\tau(\tilde{\epsilon_0}-\tilde{\epsilon})\mbox{\boldmath $V$}_e\cdot\mbox{\boldmath $R$}_e+
\tau(\mu_0-\mu)\mbox{\boldmath $V$}_m\cdot\mbox{\boldmath $R$}_m
\\
\\
\displaystyle
=\tau\tilde{\epsilon}\vert\mbox{\boldmath $R$}_e\vert^2
+\tau\mu\vert\mbox{\boldmath $R$}_m\vert^2
+\nabla\cdot(\mbox{\boldmath $R$}_e\times\mbox{\boldmath $R$}_m)
+e^{-\tau T}(\mbox{\boldmath $F$}\cdot\mbox{\boldmath $R$}_e
+\mbox{\boldmath $G$}\cdot\mbox{\boldmath $R$}_m).
\end{array}
\tag {2.10}
$$
Integrating both sides on (2.10) over $\Bbb R^3$, we obtain
$$\begin{array}{l}
\displaystyle
\,\,\,\,\,\,
\tau\int_{\Bbb R^3}
\left\{(\tilde{\epsilon_0}-\tilde{\epsilon})\mbox{\boldmath $V$}_e\cdot\mbox{\boldmath $R$}_e+
(\mu_0-\mu)\mbox{\boldmath $V$}_m\cdot\mbox{\boldmath $R$}_m\right\}dx\\
\\
\displaystyle
=\tau\int_{\Bbb R^3}(\tilde{\epsilon}\vert\mbox{\boldmath $R$}_e\vert^2
+\mu\vert\mbox{\boldmath $R$}_m\vert^2)dx
+e^{-\tau T}\int_{\Bbb R^3}(\mbox{\boldmath $F$}\cdot\mbox{\boldmath $R$}_e
+\mbox{\boldmath $G$}\cdot\mbox{\boldmath $R$}_m)dx.
\end{array}
\tag {2.11}
$$
Taking the inner product of the both sides on the first equation of (2.2) with $\mbox{\boldmath $W$}_e$,
we obtain
$$\displaystyle
\nabla\times\mbox{\boldmath $V$}_m\cdot\mbox{\boldmath $W$}_e-\tau\tilde{\epsilon_0}\mbox{\boldmath $V$}_e\cdot\mbox{\boldmath $W$}_e
+\mbox{\boldmath $f$}\cdot\mbox{\boldmath $W$}_e
=e^{-\tau T}\epsilon_0\mbox{\boldmath $E$}_0(x,T)\cdot\mbox{\boldmath $W$}_e.
\tag {2.12}
$$
Next, taking the inner product of the both sides on the first equation of (2.1) with $\mbox{\boldmath $V$}_e$,
we obtain
$$\displaystyle
\nabla\times\mbox{\boldmath $W$}_m\cdot\mbox{\boldmath $V$}_e-\tau\tilde{\epsilon}\mbox{\boldmath $W$}_e\cdot\mbox{\boldmath $V$}_e
+\mbox{\boldmath $f$}\cdot\mbox{\boldmath $V$}_e
=e^{-\tau T}\epsilon \mbox{\boldmath $E$}(x,T)\cdot\mbox{\boldmath $V$}_e.
\tag {2.13}
$$
Third, taking the inner product of the both sides on the second equation of (2.1) with $\mbox{\boldmath $V$}_m$,
we obtain
$$\displaystyle
\nabla\times\mbox{\boldmath $W$}_e\cdot\mbox{\boldmath $V$}_m+\tau\mu \mbox{\boldmath $W$}_m\cdot\mbox{\boldmath $V$}_m
=-e^{-\tau T}\mu \mbox{\boldmath $H$}(x,T)\cdot\mbox{\boldmath $V$}_m.
\tag {2.14}
$$
Fourth, taking the inner product of the both sides on the second equation of (2.2) with $\mbox{\boldmath $W$}_m$,
we obtain
$$\displaystyle
\nabla\times\mbox{\boldmath $V$}_e\cdot\mbox{\boldmath $W$}_m+\tau\mu_0 \mbox{\boldmath $V$}_m\cdot\mbox{\boldmath $W$}_m
=-e^{-\tau T}\mu_0 \mbox{\boldmath $H$}_0(x,T)\cdot\mbox{\boldmath $W$}_m.
\tag {2.15}
$$
It follows from (2.12) and (2.14) that
$$\begin{array}{l}
\displaystyle
\,\,\,\,\,\,
\nabla\cdot(\mbox{\boldmath $V$}_m\times\mbox{\boldmath $W$}_e)
-\tau\left(\tilde{\epsilon_0}\mbox{\boldmath $V$}_e\cdot\mbox{\boldmath $W$}_e
+\mu\mbox{\boldmath $W$}_m\cdot\mbox{\boldmath $V$}_m\right)
+\mbox{\boldmath $f$}\cdot\mbox{\boldmath $W$}_e\\
\\
\displaystyle
=e^{-\tau T}\left(\epsilon_0\mbox{\boldmath $E$}_0(x,T)\cdot\mbox{\boldmath $W$}_e
+\mu \mbox{\boldmath $H$}(x,T)\cdot\mbox{\boldmath $V$}_m\right)
\end{array}
$$
and thus
$$\begin{array}{ll}
\displaystyle
\int_{\Bbb R^3}\mbox{\boldmath $f$}\cdot\mbox{\boldmath $W$}_edx
&
\displaystyle
=
\tau\int_{\Bbb R^3}\left(\tilde{\epsilon_0}\mbox{\boldmath $V$}_e\cdot\mbox{\boldmath $W$}_e
+\mu\mbox{\boldmath $W$}_m\cdot\mbox{\boldmath $V$}_m\right)dx\\
\\
\displaystyle
&
\displaystyle
\,\,\,
+e^{-\tau T}
\int_{\Bbb R^3}\left(\epsilon_0\mbox{\boldmath $E$}_0(x,T)\cdot\mbox{\boldmath $W$}_e
+\mu \mbox{\boldmath $H$}(x,T)\cdot\mbox{\boldmath $V$}_m\right)dx.
\end{array}
\tag {2.16}
$$
It follows from (2.13) and (2.15) that
$$\begin{array}{l}
\displaystyle
\,\,\,\,\,\,
\nabla\cdot(\mbox{\boldmath $W$}_m\times\mbox{\boldmath $V$}_e)
-\tau\left(\tilde{\epsilon}\mbox{\boldmath $W$}_e\cdot\mbox{\boldmath $V$}_e+\mu_0\mbox{\boldmath $V$}_m\cdot\mbox{\boldmath $W$}_m\right)
+\mbox{\boldmath $f$}\cdot\mbox{\boldmath $V$}_e
\\
\\
\displaystyle
=e^{-\tau T}\left(\epsilon \mbox{\boldmath $E$}(x,T)\cdot\mbox{\boldmath $V$}_e
+\mu_0 \mbox{\boldmath $H$}_0(x,T)\cdot\mbox{\boldmath $W$}_m\right)
\end{array}
$$
and thus
$$\begin{array}{ll}
\displaystyle
\int_{\Bbb R^3}\mbox{\boldmath $f$}\cdot\mbox{\boldmath $V$}_e dx
&
\displaystyle
=\tau\int_{\Bbb R^3}\left(\tilde{\epsilon}\mbox{\boldmath $W$}_e\cdot\mbox{\boldmath $V$}_e+\mu_0\mbox{\boldmath $V$}_m\cdot\mbox{\boldmath $W$}_m\right)dx\\
\\
\displaystyle
&
\displaystyle
\,\,\,
+e^{-\tau T}\int_{\Bbb R^3}\left(\epsilon \mbox{\boldmath $E$}(x,T)\cdot\mbox{\boldmath $V$}_e
+\mu_0 \mbox{\boldmath $H$}_0(x,T)\cdot\mbox{\boldmath $W$}_m\right)dx.
\end{array}
\tag {2.17}
$$
Then from (2.11), (2.16) and (2.17) we obtain
$$\begin{array}{l}
\displaystyle
\,\,\,\,\,\,
\int_{\Bbb R^3}\mbox{\boldmath $f$}\cdot(\mbox{\boldmath $W$}_e-\mbox{\boldmath $V$}_e)dx\\
\\
\displaystyle
=\tau\int_{\Bbb R^3}\left\{(\tilde{\epsilon_0}-\tilde{\epsilon})\mbox{\boldmath $W$}_e\cdot\mbox{\boldmath $V$}_e
+(\mu-\mu_0)\mbox{\boldmath $W$}_m\cdot\mbox{\boldmath $V$}_m\right\}dx\\
\\
\displaystyle
\,\,\,
+e^{-\tau T}
\int_{\Bbb R^3}
\left(\epsilon_0\mbox{\boldmath $E$}_0(x,T)\cdot\mbox{\boldmath $W$}_e
-\epsilon\mbox{\boldmath $E$}(x,T)\cdot\mbox{\boldmath $V$}_e
+\mu \mbox{\boldmath $H$}(x,T)\cdot\mbox{\boldmath $V$}_m
-\mu_0 \mbox{\boldmath $H$}_0(x,T)\cdot\mbox{\boldmath $W$}_m\right)dx\\
\\
\displaystyle
=\tau\int_{\Bbb R^3}\left\{(\tilde{\epsilon_0}-\tilde{\epsilon})\mbox{\boldmath $V$}_e\cdot\mbox{\boldmath $V$}_e
+(\mu-\mu_0)\mbox{\boldmath $V$}_m\cdot\mbox{\boldmath $V$}_m\right\}dx\\
\\
\displaystyle
\,\,\,
+\tau\int_{\Bbb R^3}\left\{(\tilde{\epsilon_0}-\tilde{\epsilon})\mbox{\boldmath $R$}_e\cdot\mbox{\boldmath $V$}_e
+(\mu-\mu_0)\mbox{\boldmath $R$}_m\cdot\mbox{\boldmath $V$}_m\right\}dx
\\
\\
\displaystyle
\,\,\,
+e^{-\tau T}
\int_{\Bbb R^3}
\left(\epsilon_0\mbox{\boldmath $E$}_0(x,T)\cdot\mbox{\boldmath $W$}_e
-\epsilon\mbox{\boldmath $E$}(x,T)\cdot\mbox{\boldmath $V$}_e
+\mu \mbox{\boldmath $H$}(x,T)\cdot\mbox{\boldmath $V$}_m
-\mu_0 \mbox{\boldmath $H$}_0(x,T)\cdot\mbox{\boldmath $W$}_m\right)dx\\
\\
\displaystyle
=\tau\int_{\Bbb R^3}\left\{(\tilde{\epsilon_0}-\tilde{\epsilon})\vert\mbox{\boldmath $V$}_e\vert^2
+(\mu-\mu_0)\vert\mbox{\boldmath $V$}_m\vert^2\right\}dx\\
\\
\displaystyle
\,\,\,
+\tau\int_{\Bbb R^3}\left\{(\tilde{\epsilon_0}-\tilde{\epsilon})\mbox{\boldmath $V$}_e\cdot\mbox{\boldmath $R$}_e
+(\mu_0-\mu)\mbox{\boldmath $V$}_m\cdot\mbox{\boldmath $R$}_m
+2(\mu-\mu_0)\mbox{\boldmath $R$}_m\cdot\mbox{\boldmath $V$}_m\right\}dx
\\
\\
\displaystyle
\,\,\,
+e^{-\tau T}
\int_{\Bbb R^3}
\left(\epsilon_0\mbox{\boldmath $E$}_0(x,T)\cdot\mbox{\boldmath $W$}_e
-\epsilon\mbox{\boldmath $E$}(x,T)\cdot\mbox{\boldmath $V$}_e
+\mu \mbox{\boldmath $H$}(x,T)\cdot\mbox{\boldmath $V$}_m
-\mu_0 \mbox{\boldmath $H$}_0(x,T)\cdot\mbox{\boldmath $W$}_m\right)dx\\
\\
\displaystyle
=\tau\int_{\Bbb R^3}\left\{(\tilde{\epsilon_0}-\tilde{\epsilon})\vert\mbox{\boldmath $V$}_e\vert^2
+(\mu-\mu_0)\vert\mbox{\boldmath $V$}_m\vert^2\right\}dx
\\
\\
\displaystyle
\,\,\,
+\tau\int_{\Bbb R^3}(\tilde{\epsilon}\vert\mbox{\boldmath $R$}_e\vert^2
+\mu\vert\mbox{\boldmath $R$}_m\vert^2
+2(\mu-\mu_0)\mbox{\boldmath $R$}_m\cdot\mbox{\boldmath $V$}_m)dx
+e^{-\tau T}\int_{\Bbb R^3}(\mbox{\boldmath $F$}\cdot\mbox{\boldmath $R$}_e
+\mbox{\boldmath $G$}\cdot\mbox{\boldmath $R$}_m)dx
\\
\\
\displaystyle
\,\,\,
+e^{-\tau T}
\int_{\Bbb R^3}
\left(\epsilon_0\mbox{\boldmath $E$}_0(x,T)\cdot\mbox{\boldmath $W$}_e
-\epsilon\mbox{\boldmath $E$}(x,T)\cdot\mbox{\boldmath $V$}_e
+\mu \mbox{\boldmath $H$}(x,T)\cdot\mbox{\boldmath $V$}_m
-\mu_0 \mbox{\boldmath $H$}_0(x,T)\cdot\mbox{\boldmath $W$}_m\right)dx.
\end{array}
\tag {2.18}
$$
Here write
$$\displaystyle
\mu\vert\mbox{\boldmath $R$}_m\vert^2+2(\mu-\mu_0)\mbox{\boldmath $R$}_m\cdot\mbox{\boldmath $V$}_m
=\mu\left\vert\mbox{\boldmath $R$}_m+\frac{\mu-\mu_0}{\mu}\mbox{\boldmath $V$}_m\right\vert^2
-\frac{(\mu-\mu_0)^2}{\mu}\vert\mbox{\boldmath $V$}_m\vert^2.
\tag {2.19}
$$
Now from the second equation on (1.10)
and a combination of (2.6), (2.18) and (2.19) we obtain (2.7).

\noindent
$\Box$

\proclaim{\noindent Lemma 2.1.}
We have
$$\displaystyle
\Vert\mbox{\boldmath $V$}_e\Vert_{L^2(\Bbb R^3)}
=\Vert\mbox{\boldmath $V$}_m\Vert_{L^2(\Bbb R^3)}
=O(\tau^{-5/2}).
\tag {2.20}
$$
and
$$\displaystyle
\Vert\mbox{\boldmath $W$}_e\Vert_{L^2(\Bbb R^3)}
=\Vert\mbox{\boldmath $W$}_m\Vert_{L^2(\Bbb R^3)}
=O(\tau^{-5/2}).
\tag {2.21}
$$

\endproclaim

{\it\noindent Proof.}
Taking the inner product of the first equation on (2.2) with $\mbox{\boldmath $V$}_e$ we obtain
$$\displaystyle
\nabla\times\mbox{\boldmath $V$}_m\cdot\mbox{\boldmath $V$}_e-\tau\tilde{\epsilon_0}\vert\mbox{\boldmath $V$}_e\vert^2+
\mbox{\boldmath $f$}\cdot\mbox{\boldmath $V$}_e
=e^{-\tau T}\epsilon_0\mbox{\boldmath $E$}_0(x,T)\cdot\mbox{\boldmath $V$}_e.
\tag {2.22}
$$
Taking the inner product of the second equation on (2.2) with $\mbox{\boldmath $V$}_m$ we obtain
$$\displaystyle
\nabla\times\mbox{\boldmath $V$}_e\cdot\mbox{\boldmath $V$}_m+\tau\mu_0\vert\mbox{\boldmath $V$}_m\vert^2
=-e^{-\tau T}\mu_0\mbox{\boldmath $H$}_0(x,T)\cdot\mbox{\boldmath $V$}_m.
\tag {2.23}
$$
A combination of (2.22) and (2.23) gives
$$\displaystyle
\nabla\cdot(\mbox{\boldmath $V$}_m\times\mbox{\boldmath $V$}_e)
-\tau\left(\tilde{\epsilon_0}\vert\mbox{\boldmath $V$}_e\vert^2
+\mu_0\vert\mbox{\boldmath $V$}_m\vert^2\right)
+\mbox{\boldmath $f$}\cdot\mbox{\boldmath $V$}_e
=e^{-\tau T}
\left(\epsilon_0\mbox{\boldmath $E$}_0(x,T)\cdot\mbox{\boldmath $V$}_e
+\mu_0\mbox{\boldmath $H$}_0(x,T)\cdot\mbox{\boldmath $V$}_m\right).
$$
Thus, integrating both sides over $\Bbb R^3$ we obtain
$$\begin{array}{l}
\displaystyle
\,\,\,\,\,\,
\tau\int_{\Bbb R^3}\left(\tilde{\epsilon_0}\vert\mbox{\boldmath $V$}_e\vert^2
+\mu_0\vert\mbox{\boldmath $V$}_m\vert^2\right)dx\\
\\
\displaystyle
\,\,\,
+\int_{\Bbb R^3}\left\{(e^{-\tau T}\epsilon_0\mbox{\boldmath $E$}_0(x,T)-\mbox{\boldmath $f$})\cdot\mbox{\boldmath $V$}_e
+e^{-\tau T}\mu_0\mbox{\boldmath $H$}_0(x,T)\cdot\mbox{\boldmath $V$}_m\right\}dx
\\
\\
\displaystyle
=0,
\end{array}
$$
that is
$$\begin{array}{l}
\displaystyle
\,\,\,\,\,\,
\int_{\Bbb R^3}
\left(\tilde{\epsilon_0}
\left\vert\mbox{\boldmath $V$}_e+
\frac{e^{-\tau T}\epsilon_0\mbox{\boldmath $E$}_0(x,T)-\mbox{\boldmath $f$}}{2\tau\tilde{\epsilon_0}}\right\vert^2
+\mu_0
\left\vert\mbox{\boldmath $V$}_m+
\frac{e^{-\tau T}\mbox{\boldmath $H$}_0(x,T)}{2\tau}\right\vert^2\right)
dx\\
\\
\displaystyle
=\frac{1}{4\tau^2}\int_{\Bbb R^3}\left(\tilde{\epsilon_0}^{-1}\vert e^{-\tau T}\epsilon_0\mbox{\boldmath $E$}_0-\mbox{\boldmath $f$}\vert^2
+\mu_0e^{-2\tau T}\vert\mbox{\boldmath $H$}_0(x,T)\vert^2\right)dx.
\end{array}
\tag {2.24}
$$
Since
$$\displaystyle
\left\vert\mbox{\boldmath $V$}_e+
\frac{e^{-\tau T}\epsilon_0\mbox{\boldmath $E$}_0(x,T)-\mbox{\boldmath $f$}}{2\tau\tilde{\epsilon_0}}\right\vert^2
\ge \frac{1}{2}\vert\mbox{\boldmath $V$}_e\vert^2-
\frac{1}{4\tau^2\tilde{\epsilon_0}^2}\vert e^{-\tau T}\epsilon_0\mbox{\boldmath $E$}_0(x,T)-\mbox{\boldmath $f$}\vert^2
$$
and
$$\displaystyle
\left\vert\mbox{\boldmath $V$}_m+
\frac{e^{-\tau T}\mbox{\boldmath $H$}_0(x,T)}{2\tau}\right\vert^2
\ge\frac{1}{2}\vert\mbox{\boldmath $V$}_m\vert^2
-\frac{1}{4\tau^2}e^{-2\tau T}\vert\mbox{\boldmath $H$}_0(x,T)\vert^2,
$$
from (2.24) we obtain
$$\begin{array}{l}
\displaystyle
\,\,\,\,\,\,
\frac{1}{2}\int_{\Bbb R^3}
\left(\tilde{\epsilon_0}\vert\mbox{\boldmath $V$}_e\vert^2
+\mu_0\vert\mbox{\boldmath $V$}_m\vert^2\right)dx
\\
\\
\displaystyle
\le\frac{1}{2\tau^2}\int_{\Bbb R^3}\left(\tilde{\epsilon_0}^{-1}\vert e^{-\tau T}\epsilon_0\mbox{\boldmath $E$}_0-\mbox{\boldmath $f$}\vert^2
+\mu_0e^{-2\tau T}\vert\mbox{\boldmath $H$}_0(x,T)\vert^2\right)dx,
\end{array}
$$
that is
$$\begin{array}{l}
\displaystyle
\,\,\,\,\,\,
\int_{\Bbb R^3}
\left(\tilde{\epsilon_0}\vert\mbox{\boldmath $V$}_e\vert^2
+\mu_0\vert\mbox{\boldmath $V$}_m\vert^2\right)dx
\\
\\
\displaystyle
\le\frac{1}{\tau^2}\int_{\Bbb R^3}\left(\tilde{\epsilon_0}^{-1}\vert e^{-\tau T}\epsilon_0\mbox{\boldmath $E$}_0-\mbox{\boldmath $f$}\vert^2
+\mu_0e^{-2\tau T}\vert\mbox{\boldmath $H$}_0(x,T)\vert^2\right)dx.
\end{array}
$$
Therefore we obtain, as $\tau\longrightarrow\infty$
$$
\displaystyle
\int_{\Bbb R^3}
\left(\tilde{\epsilon_0}\vert\mbox{\boldmath $V$}_e\vert^2
+\mu_0\vert\mbox{\boldmath $V$}_m\vert^2\right)dx
=O(\tau^{-2}(e^{-2\tau T}+\Vert\mbox{\boldmath $f$}\Vert_{L^2(\Bbb R^3)}^2))
$$
and thus from (1.5) we obtain (2.20).

Since $\mbox{\boldmath $W$}_e$ and $\mbox{\boldmath $W$}_m$ satisfy (2.1),
similarly to the derivation of (2.20), we obtain
$$\displaystyle
\Vert\mbox{\boldmath $W$}_e\Vert_{L^2(\Bbb R^3)}
=\Vert\mbox{\boldmath $W$}_m\Vert_{L^2(\Bbb R^3)}
=O(\tau^{-1}(e^{-\tau T}+\Vert\mbox{\boldmath $f$}\Vert_{L^2(\Bbb R^3)})).
$$
This together with (1.5) yields (2.21).

\noindent
$\Box$

\proclaim{\noindent Lemma 2.2.}
We have, as $\tau\longrightarrow\infty$
$$\begin{array}{l}
\displaystyle
\,\,\,\,\,\,
\int_{\Bbb R^3}\left(\tilde{\epsilon}\vert\mbox{\boldmath $R$}_e\vert^2+
\mu\left\vert\mbox{\boldmath $R$}_m+\frac{\mu-\mu_0}{\mu}\mbox{\boldmath $V$}_m\right\vert^2\right)dx
\\
\\
\displaystyle
\le 
\int_{\Bbb R^3}\left\{\frac{(\tilde{\epsilon_0}-\tilde{\epsilon})^2}{\tilde{\epsilon}}\vert\mbox{\boldmath $V$}_e\vert^2+
\frac{(\mu-\mu_0)^2}{\mu}\vert\mbox{\boldmath $V$}_m\vert^2\right\}dx
+O(\tau^{-7/2}e^{-\tau T}).
\end{array}
\tag {2.25}
$$

\endproclaim
{\it\noindent Proof.}
It follows from (2.11) that
$$\begin{array}{l}
\displaystyle
\int_{\Bbb R^3}\tilde{\epsilon}\vert\mbox{\boldmath $R$}_e\vert^2
+\left(\frac{e^{-\tau T}}{\tau}\mbox{\boldmath $F$}-(\tilde{\epsilon_0}-\tilde{\epsilon})
\mbox{\boldmath $V$}e\right)\cdot\mbox{\boldmath $R$}_e\,dx
\\
\\
\displaystyle
\,\,\,
+\int_{\Bbb R^3}
\mu\vert\mbox{\boldmath $R$}_m\vert^2+
\left((\mu-\mu_0)\mbox{\boldmath $V$}_m
+\frac{e^{-\tau T}}{\tau}\mbox{\boldmath $G$}\right)\cdot\mbox{\boldmath $R$}_m\,dx
\\
\\
\displaystyle
=0.
\end{array}
$$

Rewrite this as
$$\begin{array}{l}
\displaystyle
\,\,\,\,\,\,
\int_{\Bbb R^3}
\tilde{\epsilon}\left\vert\mbox{\boldmath $R$}_e+\frac{1}{2\tilde{\epsilon}}
\left(\frac{e^{-\tau T}\mbox{\boldmath $F$}}{\tau}-(\tilde{\epsilon_0}-\tilde{\epsilon})\mbox{\boldmath $V$}_e\right)\right\vert^2\,dx
\\
\\
\displaystyle
\,\,\,
+\int_{\Bbb R^3}
\mu\left\vert\mbox{\boldmath $R$}_m+\frac{\mu-\mu_0}{2\mu}\mbox{\boldmath $V$}_m+
\frac{e^{-\tau T}\mbox{\boldmath $G$}}{2\mu\tau}\right\vert^2\,
dx\\
\\
\displaystyle
=\frac{1}{4}
\int_{\Bbb R^3}
\left\{
\frac{1}{\tilde{\epsilon}}
\left\vert\frac{e^{-\tau T}\mbox{\boldmath $F$}}{\tau}-(\tilde{\epsilon_0}-\tilde{\epsilon})\mbox{\boldmath $V$}_e\right\vert^2
+\frac{1}{\mu}
\left\vert
(\mu-\mu_0)\mbox{\boldmath $V$}_m
+\frac{e^{-\tau T}\mbox{\boldmath $G$}}{\tau}\right\vert^2
\right\}dx.
\end{array}
$$
and
$$\displaystyle
\left\vert\mbox{\boldmath $R$}_m+\frac{\mu-\mu_0}{2\mu}\mbox{\boldmath $V$}_m+
\frac{e^{-\tau T}\mbox{\boldmath $G$}}{2\mu\tau}\right\vert^2
=\left\vert\left(\mbox{\boldmath $R$}_m+\frac{\mu-\mu_0}{\mu}\mbox{\boldmath $V$}_m\right)
-\frac{1}{2\mu}
\left((\mu-\mu_0)\mbox{\boldmath $V$}_m-
\frac{e^{-\tau T}\mbox{\boldmath $G$}}{\tau}
\right)\right\vert^2.
$$
Applying the inequality $\vert a+b\vert^2\ge\frac{1}{2}\vert a\vert^2-\vert b\vert^2$ to the left-hand side and using (2.20), we obtain
$$\begin{array}{l}
\displaystyle
\,\,\,\,\,\,
\frac{1}{2}\int_{\Bbb R^3}\left(\tilde{\epsilon}\vert\mbox{\boldmath $R$}_e\vert^2+
\mu\left\vert\mbox{\boldmath $R$}_m+\frac{\mu-\mu_0}{\mu}\mbox{\boldmath $V$}_m\right\vert^2\right)dx
\\
\\
\displaystyle
\le\frac{1}{2}
\int_{\Bbb R^3}
\frac{1}{\tilde{\epsilon}}
\left\vert\frac{e^{-\tau T}\mbox{\boldmath $F$}}{\tau}-(\tilde{\epsilon_0}-\tilde{\epsilon})\mbox{\boldmath $V$}_e\right\vert^2
\,dx
\\
\\
\displaystyle
\,\,\,
+\frac{1}{4}
\int_{\Bbb R^3}
\frac{1}{\mu}
\left\{
\left\vert
(\mu-\mu_0)\mbox{\boldmath $V$}_m+
\frac{e^{-\tau T}\mbox{\boldmath $G$}}{\tau}\right\vert^2
+\left\vert(\mu-\mu_0)\mbox{\boldmath $V$}_m-\frac{e^{-\tau T}\mbox{\boldmath $G$}}{\tau}\right\vert^2
\right\}\,dx
\\
\\
\displaystyle
=O(\tau^{-2}e^{-2\tau T})+O(\tau^{-1}e^{-\tau T}\tau^{-5/2})+
\int_{\Bbb R^3}\,\left\{\frac{(\tilde{\epsilon_0}-\tilde{\epsilon})^2}{2\tilde{\epsilon}}\vert\mbox{\boldmath $V$}_e\vert^2
+\frac{(\mu-\mu_0)^2}{2\mu}\vert\mbox{\boldmath $V$}_m\vert^2\right\}\,dx.
\end{array}
$$
Therefore, as $\tau\longrightarrow\infty$ we obtain (2.25).

\noindent
$\Box$

Now we are ready to prove Theorem 1.1.

Applying (2.20), (2.21) and (2.25) to the right-hand side on (2.7)
and the equations (1.10)
we obtain (1.6).  The estimate (1.7) is clearly valid.  This completes the proof of Theorem 1.1.

\section{Proof of Theorem 1.2}

\subsection{A reduction to time independent case}

In this subsection, we do not assume that $\epsilon_0$, $\mu_0$ and $\sigma_0$ are constant.

Let $\mbox{\boldmath $V$}_e^0$ and $\mbox{\boldmath $V$}_m^0$ solve
$$
\left\{
\begin{array}{ll}
\displaystyle
\nabla\times\mbox{\boldmath $V$}_m^0-\tau\tilde{\epsilon_0}\mbox{\boldmath $V$}_e^0+\mbox{\boldmath $f$}
=\mbox{\boldmath $0$},
&
\displaystyle
x\in\Bbb R^3,
\\
\\
\displaystyle
\nabla\times\mbox{\boldmath $V$}_e^0+\tau\mu_0\mbox{\boldmath $V$}_m^0
=\mbox{\boldmath $0$},
&
\displaystyle
x\in\Bbb R^3.
\end{array}
\right.
\tag {3.1}
$$
Using a similar argument to derive (2.20) and (2.25), we have
$$\displaystyle
\Vert\mbox{\boldmath $V$}_e^0\Vert_{L^2(\Bbb R^3)}
=\Vert\mbox{\boldmath $V$}_m^0\Vert_{L^2(\Bbb R^3)}
=O(\tau^{-5/2}).
\tag {3.2}
$$
and
$$
\displaystyle
\,\,\,\,\,\,
\int_{\Bbb R^3}\left(\tilde{\epsilon_0}\vert\mbox{\boldmath $V$}_e-\mbox{\boldmath $V$}_e^0\vert^2+
\mu_0\vert\mbox{\boldmath $V$}_m-\mbox{\boldmath $V$}_m^0\vert^2\right)dx
=O(\tau^{-2}e^{-2\tau T}).
\tag {3.3}
$$
Writing
$$\left\{
\begin{array}{ll}
\displaystyle
\mbox{\boldmath $V$}_e=\mbox{\boldmath $V$}_e^0+(\mbox{\boldmath $V$}_e-\mbox{\boldmath $V$}_e^0),\\
\\
\displaystyle
\mbox{\boldmath $V$}_m=\mbox{\boldmath $V$}_m^0+(\mbox{\boldmath $V$}_m-\mbox{\boldmath $V$}_m^0),
\end{array}
\right.
$$
and applying (3.2) and (3.3), from Theorem 1.1 we obtain the following asymptotic estimates.

\proclaim{\noindent Proposition 3.1.}
We have, as $\tau\longrightarrow\infty$
$$
\displaystyle
I(\tau)
\le 
\tau\int_{\Bbb R^3}
\left\{
\frac{\tilde{\epsilon_0}}{\tilde{\epsilon}}(\tilde{\epsilon_0}-\tilde{\epsilon})
\vert\mbox{\boldmath $V$}_e^0\vert^2dx
+
(\mu-\mu_0)
\vert\mbox{\boldmath $V$}_m^0\vert^2
\right\}dx
+O(\tau^{-5/2}e^{-\tau T})
\tag {3.4}
$$
and
$$\displaystyle
I(\tau)
\ge 
\tau\int_{\Bbb R^3}\left\{(\tilde{\epsilon_0}-\tilde{\epsilon})\vert\mbox{\boldmath $V$}_e^0\vert^2
+\frac{\mu_0}{\mu}(\mu-\mu_0)\vert\mbox{\boldmath $V$}_m^0\vert^2\right\}dx
+O(\tau^{-5/2}e^{-\tau T}).
\tag {3.5}
$$

\endproclaim

{\bf\noindent Remark 3.1.}
It follows from (3.3) that
$$
\displaystyle
\Vert\mbox{\boldmath $V$}_e-\mbox{\boldmath $V$}_e^0\Vert_{L^2(\Bbb R^3)}
=O(\tau^{-1}e^{-\tau T})
$$
and this together with (1.5) yields
$$\displaystyle
I(\tau)
=\tilde{I}(\tau)
+O(\tau^{-5/2}e^{-\tau T}),
$$
where
$$\begin{array}{ll}
\displaystyle
\tilde{I}(\tau)
=
\int_{\Bbb R^3}\mbox{\boldmath $f$}\cdot(\mbox{\boldmath $W$}_e-\mbox{\boldmath $V$}_e^0)dx,
& \displaystyle
\tau>0.
\end{array}
$$
Therefore from (3.4) and (3.5) we have, as $\tau\longrightarrow\infty$
$$
\displaystyle
\tilde{I}(\tau)
\le 
\tau\int_{\Bbb R^3}
\left\{
\frac{\tilde{\epsilon_0}}{\tilde{\epsilon}}(\tilde{\epsilon_0}-\tilde{\epsilon})
\vert\mbox{\boldmath $V$}_e^0\vert^2dx
+
(\mu-\mu_0)
\vert\mbox{\boldmath $V$}_m^0\vert^2
\right\}dx
+O(\tau^{-5/2}e^{-\tau T})
$$
and
$$\displaystyle
\tilde{I}(\tau)
\ge 
\tau\int_{\Bbb R^3}\left\{(\tilde{\epsilon_0}-\tilde{\epsilon})\vert\mbox{\boldmath $V$}_e^0\vert^2
+\frac{\mu_0}{\mu}(\mu-\mu_0)\vert\mbox{\boldmath $V$}_m^0\vert^2\right\}dx
+O(\tau^{-5/2}e^{-\tau T}).
$$
This means that, instead of the original indicator function, it is possible to use
another indicator function $\tilde{I}(\tau)$.

Needless to say, when $\epsilon_0$, $\mu_o$ and $\sigma_0$ are constant, one obtains the completely same result as Theorem 1.2
for this indicator function.

\subsection{The solutions of the system (3.1) in the case when $\epsilon_0$, $\mu_0$ and $\sigma_0$ are constant}

In this subsection we assume that $\epsilon_0$, $\mu_0$ and $\sigma_0$ are constant.
Then, the system (3.1) is equivalent to the equation
$$\displaystyle
\nabla\times\nabla\times\mbox{\boldmath $V$}_e^0+\tau^2\tilde{\epsilon_0}\mu_0\mbox{\boldmath $V$}_e^0-\tau\mu_0\mbox{\boldmath $f$}
=\mbox{\boldmath $0$}
\tag {3.6}
$$
and
$$\displaystyle
\mbox{\boldmath $V$}_m^0=-\frac{1}{\tau\mu_0}\nabla\times\mbox{\boldmath $V$}_e^0.
\tag {3.7}
$$
Let us recall how to construct the solutions of equation (3.6), which has been done in \cite{IMax}.

Hereafter we simply write $\mbox{\boldmath $V$}_e^0=\mbox{\boldmath $V$}$.
Following \cite{ABF}, we write the solution of (3.6) in the form
$$\displaystyle
\mbox{\boldmath $V$}=\mbox{\boldmath $V$}_0+\mbox{\boldmath $V$}_1
$$
Then, one has the expression
$$\begin{array}{l}
\displaystyle
\,\,\,\,\,\,\,
\nabla\times\nabla\times\mbox{\boldmath $V$}+\tau^2\tilde{\epsilon_0}\mu_0
\mbox{\boldmath $V$}-\tau\mu_0\mbox{\boldmath $f$}(x,\tau)\\
\\
\displaystyle
=\left\{-(\Delta-\tau^2\tilde{\epsilon_0}\mu_0)\mbox{\boldmath $V$}_0-\tau\mu_0\mbox{\boldmath $f$}(x,\tau)\right\}
+\left(\tau^2\tilde{\epsilon_0}\mu_0\mbox{\boldmath $V$}_1+\nabla(\nabla\cdot\mbox{\boldmath $V$}_0)\right)
+\nabla\times\nabla\times\mbox{\boldmath $V$}_1.
\end{array}
$$
From this we see that if
$$\displaystyle
(\Delta-\tau^2\tilde{\epsilon_0}\mu_0)\mbox{\boldmath $V$}_0+\tau\mu_0\mbox{\boldmath $f$}(x,\tau)=\mbox{\boldmath $0$}
\tag {3.8}
$$
and
$$\displaystyle
\tau^2\tilde{\epsilon_0}\mu_0\mbox{\boldmath $V$}_1+\nabla(\nabla\cdot\mbox{\boldmath $V$}_0)=\mbox{\boldmath $0$},
\tag {3.9}
$$
then $\nabla\times\mbox{\boldmath $V$}_1=\mbox{\boldmath $0$}$, and thus $\mbox{\boldmath $V$}=\mbox{\boldmath $V$}_0+\mbox{\boldmath $V$}_1$
satisfies (3.6).  

Note also from (3.9) we have
$$\displaystyle
\mbox{\boldmath $V$}_e^0=\mbox{\boldmath $V$}_0-(\tau^2\tilde{\epsilon_0}\mu_0)^{-1}\nabla(\nabla\cdot\mbox{\boldmath $V$}_0)
\tag {3.10}
$$
and the explicit form of $\mbox{\boldmath $V$}_0\in L^2(\Bbb R^3)^3$ satisfying (3.8).

By (18) in \cite{IMax} we know that $\mbox{\boldmath $V$}_e^0$ given by (3.10) is smooth outside $B$ and has the form
$$\begin{array}{ll}
\displaystyle
\mbox{\boldmath $V$}_e^0
=K(\tau)\tilde{f}(\tau)v\mbox{\boldmath $M$}\mbox{\boldmath $a$}
& x\in\Bbb R^3\setminus\overline B,
\end{array}
\tag {3.11}
$$
where $v=v(x,\tau)$ has the form
$$\displaystyle
v(x,\tau)=\frac{e^{-\tau\sqrt{\mu_0\tilde{\epsilon_0}}\,\vert x-p\vert}}{\vert x-p\vert},
$$
$$
\left\{
\begin{array}{l}
\displaystyle
K(\tau)=\frac{\mu_0\tau\varphi(\tau\sqrt{\mu_0\tilde{\epsilon_0}}\,\eta)}{(\tau\sqrt{\mu_0\tilde{\epsilon_0}})^3},
\\
\\
\displaystyle
\varphi(\xi)=\xi\cosh\xi-\sinh\xi
\end{array}
\right.
$$
and
$$
\left\{
\begin{array}{l}
\displaystyle
\mbox{\boldmath $M$}
=\mbox{\boldmath $M$}(x;\tau)
=AI_3-B\,\frac{x-p}{\vert x-p\vert}\otimes
\frac{x-p}{\vert x-p\vert},
\\
\\
\displaystyle
A=A(x,\tau)=1+\frac{1}{\tau\sqrt{\mu_0\tilde{\epsilon_0}}}
\left(\frac{1}{\vert x-p\vert}
+\frac{1}{\tau\sqrt{\mu_0\tilde{\epsilon_0}}\vert x-p\vert^2}\right),\\
\\
\displaystyle
B=B(x,\tau)=1+\frac{3}{\tau\sqrt{\mu_0\tilde{\epsilon_0}}}
\left(\frac{1}{\vert x-p\vert}
+\frac{1}{\tau\sqrt{\mu_0\tilde{\epsilon_0}}\vert x-p\vert^2}\right).
\end{array}
\right.
$$
And from (3.7) and (3.11) we obtain
$$\begin{array}{ll}
\displaystyle
\mbox{\boldmath $V$}_m^0=-\frac{1}{\tau\mu_0}K(\tau)\tilde{f}(\tau)\nabla v\times(\mbox{\boldmath $M$}\mbox{\boldmath $a$}),
&
x\in\Bbb R^3\setminus\overline B.
\end{array}
\tag {3.12}
$$
The expression (3.11) is a simple application of the mean value theorem \cite{CH} for the modified
Helmholtz equation to the explicit form of $\mbox{\boldmath $V$}_0$.

\proclaim{\noindent Lemma 3.1.}  Let $x\in\Bbb R^3\setminus\overline B$.
We have
$$\displaystyle
\vert\mbox{\boldmath $V$}_e^0\vert^2
=
K(\tau)^2(\tilde{f}(\tau))^2v^2
\left\{A^2\left\vert\mbox{\boldmath $a$}\times\frac{x-p}{\vert x-p\vert}\right\vert^2
+(B-A)^2\left(\mbox{\boldmath $a$}\cdot\frac{x-p}{\vert x-p\vert}\right)^2\right\}
\tag {3.13}
$$
and
$$\displaystyle
\vert\mbox{\boldmath $V$}_m^0\vert^2
=K(\tau)^2\tilde{f}(\tau)^2v^2\left(\sqrt{\frac{\tilde{\epsilon_0}}{\mu_0}}+\frac{1}{\tau\mu_0\vert x-p\vert}\right)^2A^2\left\vert\mbox{\boldmath $a$}\times\frac{x-p}{\vert x-p\vert}\right\vert^2.
\tag {3.14}
$$

\endproclaim

{\it\noindent Proof.}
Since $M^T=M$ and 
$$\displaystyle
M^2=A^2I_3+(B^2-2AB)\frac{x-p}{\vert x-p\vert}\otimes
\frac{x-p}{\vert x-p\vert},
$$
we have
$$\begin{array}{ll}
\displaystyle
\vert M\mbox{\boldmath $a$}\vert^2
&
\displaystyle
=M^tM\mbox{\boldmath $a$}\cdot\mbox{\boldmath $a$}
\\
\\
\displaystyle
&
\displaystyle
=A^2+(B^2-2AB)\left(\mbox{\boldmath $a$}\cdot\frac{x-p}{\vert x-p\vert}\right)^2
\\
\\
\displaystyle
&
\displaystyle
=A^2\left\vert\mbox{\boldmath $a$}\times\frac{x-p}{\vert x-p\vert}\right\vert^2
+(B-A)^2\left(\mbox{\boldmath $a$}\cdot\frac{x-p}{\vert x-p\vert}\right)^2.
\end{array}
$$
Thus, from (3.11) one gets (3.13).

On the other hand, it follows from (3.12) that
$$\begin{array}{l}
\displaystyle
\vert\mbox{\boldmath $V$}_m^0\vert^2
=(\tau\mu_0)^{-2}K(\tau)^2\tilde{f}(\tau)^2\vert\nabla v\times(\mbox{\boldmath $M$}\mbox{\boldmath $a$})\vert^2.
\end{array}
$$

Here we have
$$\displaystyle
\vert\nabla v\times(\mbox{\boldmath $M$}\mbox{\boldmath $a$})\vert^2
=\vert\nabla v\vert^2\vert\mbox{\boldmath $M$}\mbox{\boldmath $a$}\vert^2
-(\mbox{\boldmath $M$}\nabla v\cdot\mbox{\boldmath $a$})^2,
$$
$$\displaystyle
\nabla v=-\left(\tau\sqrt{\mu_0\tilde{\epsilon_0}}+\frac{1}{\vert x-p\vert}\right)
\frac{x-p}{\vert x-p\vert}v
$$
and
$$\displaystyle
\mbox{\boldmath $M$}\frac{x-p}{\vert x-p\vert}
=(A-B)\frac{x-p}{\vert x-p\vert}.
$$

Thus, we have
$$\displaystyle
\mbox{\boldmath $M$}\nabla v=-(A-B)\left(\tau\sqrt{\mu_0\tilde{\epsilon_0}}+\frac{1}{\vert x-p\vert}\right)
\frac{x-p}{\vert x-p\vert}v
$$
and hence
$$\displaystyle
(\mbox{\boldmath $M$}\nabla v\cdot\mbox{\boldmath $a$})^2
=(A-B)^2\left(\tau\sqrt{\mu_0\tilde{\epsilon_0}}+\frac{1}{\vert x-p\vert}\right)^2
\left(\mbox{\boldmath $a$}\cdot\frac{x-p}{\vert x-p\vert}\right)^2v^2.
$$

Thus one gets
$$
\displaystyle
\vert\nabla v\times(\mbox{\boldmath $M$}\mbox{\boldmath $a$})\vert^2
=\left(\tau\sqrt{\mu_0\tilde{\epsilon_0}}+\frac{1}{\vert x-p\vert}\right)^2v^2A^2\left\vert\mbox{\boldmath $a$}\times\frac{x-p}{\vert x-p\vert}\right\vert^2.
$$

Therefore we have the expression
$$\displaystyle
\vert\mbox{\boldmath $V$}_m^0\vert^2
=(\tau\mu_0)^{-2}K(\tau)^2\tilde{f}(\tau)^2v^2\left(\tau\sqrt{\mu_0\tilde{\epsilon_0}}+\frac{1}{\vert x-p\vert}\right)^2A^2\left\vert\mbox{\boldmath $a$}\times\frac{x-p}{\vert x-p\vert}\right\vert^2.
$$
This yields (3.14).

\noindent
$\Box$

\subsection{Finishing the proof of Theorem 1.2}

From (3.11) and (3.12) together with Proposition 3.1, one immediately see that the statement (i) in Theorem 1.2 is valid.

The following two lemmata clarify the meaning of the jump conditions (A.I) and (A.II) in the statement (ii).
Once we have those lemmata, using (1.13) and Proposition 3.1 we immediately finish the proof of (ii).

\proclaim{\noindent Lemma 3.2.}
(i)  If there exists a positive constant $C_4$ such that, for almost all $x\in D$
$$\displaystyle
\frac{\epsilon_0}{\epsilon}(\epsilon_0-\epsilon)+
(\mu-\mu_0)\cdot\frac{\epsilon_0}{\mu_0}\le -C_4,
\tag {3.15}
$$
then, we have
$$\begin{array}{l}
\displaystyle
\,\,\,\,\,\,
\int_{\Bbb R^3}
\left\{\frac{\tilde{\epsilon_0}}{\tilde{\epsilon}}(\tilde{\epsilon_0}-\tilde{\epsilon})
\vert\mbox{\boldmath $V$}_e^0\vert^2
+
(\mu-\mu_0)
\vert\mbox{\boldmath $V$}_m^0\vert^2
\right\}\,dx
\\
\\
\displaystyle
\le -C_5\tau^{-\kappa}K(\tau)^2(\tilde{f}(\tau))^2e^{-2\tau\sqrt{\mu_0\epsilon_0}\,\text{dist}\,(D,B)},
\end{array}
\tag {3.16}
$$
where $C_5$ is a positive constant and
$$\displaystyle
\kappa=
\left\{
\begin{array}{ll}
\displaystyle
3 & \text{if (B.I) or (B.II) is satisfied,}
\\
\\
\displaystyle
2 & \text{if (B.III) is satisfied.}
\end{array}
\right.
$$

(ii)  If there exists a positive constant $C_6$ such that, for almost all $x\in D$
$$\displaystyle
(\epsilon_0-\epsilon)+
\frac{\epsilon_0}{\mu}(\mu-\mu_0)
\ge C_6,
\tag {3.17}
$$
then, we have
$$
\displaystyle
\int_{\Bbb R^3}
\left\{
(\tilde{\epsilon_0}-\tilde{\epsilon})\vert\mbox{\boldmath $V$}_e^0\vert^2+
\frac{\mu_0}{\mu}(\mu-\mu_0)\vert\mbox{\boldmath $V$}_m^0\vert^2
\right\}\,dx
\ge C_7\tau^{-\kappa}K(\tau)^2(\tilde{f}(\tau))^2e^{-2\tau\sqrt{\mu_0\epsilon_0}\,\text{dist}\,(D,B)},
\tag {3.18}
$$
where $C_7$ is a positive constant and $\kappa$ is the same constant as (i).

\endproclaim

\proclaim{\noindent Lemma 3.3.}
(3.15) and (3.17) are equivalent to (A.I) and (A.II), respectively.

\endproclaim

{\it\noindent Proof of Lemma 3.2.}
It follows from (3.13) and (3.14) that
$$\begin{array}{l}
\,\,\,\,\,\,
\displaystyle
\frac{\tilde{\epsilon_0}}{\tilde{\epsilon}}(\tilde{\epsilon_0}-\tilde{\epsilon})
\vert\mbox{\boldmath $V$}_e^0\vert^2
+
(\mu-\mu_0)
\vert\mbox{\boldmath $V$}_m^0\vert^2\\
\\
\displaystyle
=\left\{
\frac{\tilde{\epsilon_0}}{\tilde{\epsilon}}(\tilde{\epsilon_0}-\tilde{\epsilon})+
(\mu-\mu_0)
\left(\sqrt{\frac{\tilde{\epsilon}_0}{\mu_0}}+\frac{1}{\tau\mu_0\vert x-p\vert}\right)^2
\right\}K(\tau)^2(\tilde{f}(\tau))^2v^2A^2
\left\vert\mbox{\boldmath $a$}\times\frac{x-p}{\vert x-p\vert}\right\vert^2
\\
\\
\displaystyle
\,\,\,
+\frac{\tilde{\epsilon_0}}{\tilde{\epsilon}}(\tilde{\epsilon_0}-\tilde{\epsilon})
K(\tau)^2(\tilde{f}(\tau)^2v^2(B-A)^2\left(\mbox{\boldmath $a$}\cdot\frac{x-p}{\vert x-p\vert}\right)^2.
\end{array}
\tag {3.19}
$$
Here from (1.8) and (1.9) we have
$$\displaystyle
\frac{\tilde{\epsilon_0}}{\tilde{\epsilon}}(\tilde{\epsilon_0}-\tilde{\epsilon})
=\frac{\epsilon_0}{\epsilon}(\epsilon_0-\epsilon)+O(\frac{1}{\tau}),
$$
$$\begin{array}{l}
\displaystyle
\,\,\,\,\,\,
\frac{\tilde{\epsilon_0}}{\tilde{\epsilon}}(\tilde{\epsilon_0}-\tilde{\epsilon})+
(\mu-\mu_0)
\left(\sqrt{\frac{\tilde{\epsilon}_0}{\mu_0}}+\frac{1}{\tau\mu_0\vert x-p\vert}\right)^2
\\
\\
\displaystyle
=\frac{\epsilon_0}{\epsilon}(\epsilon_0-\epsilon)+
(\mu-\mu_0)\cdot\frac{\epsilon_0}{\mu_0}
+O(\frac{1}{\tau})
\end{array}
$$
and
$$\left\{
\begin{array}{l}
\displaystyle
A^2=1+O(\frac{1}{\tau}),
\\
\\
\displaystyle
(A-B)^2=\frac{4}{\tau^2\mu_0\epsilon_0}\cdot\frac{1}{\vert x-p\vert^2}\left(1+O(\frac{1}{\tau})\right),
\end{array}
\right.
$$
uniformly with repect to $x\in D$.

Thus, from (3.19) one gets
$$\begin{array}{l}
\,\,\,\,\,\,
\displaystyle
\frac{\tilde{\epsilon_0}}{\tilde{\epsilon}}(\tilde{\epsilon_0}-\tilde{\epsilon})
\vert\mbox{\boldmath $V$}_e^0\vert^2
+
(\mu-\mu_0)
\vert\mbox{\boldmath $V$}_m^0\vert^2\\
\\
\displaystyle
=\left\{\frac{\epsilon_0}{\epsilon}(\epsilon_0-\epsilon)+
(\mu-\mu_0)\cdot\frac{\epsilon_0}{\mu_0}
+O(\frac{1}{\tau})\right\}K(\tau)^2(\tilde{f}(\tau))^2v^2
\left\vert\mbox{\boldmath $a$}\times\frac{x-p}{\vert x-p\vert}\right\vert^2
\\
\\
\displaystyle
\,\,\
+
\frac{4}{\tau^2\mu_0\epsilon_0\vert x-p\vert^2}
K(\tau)^2(\tilde{f}(\tau))^2v^2
\left(\frac{\epsilon_0}{\epsilon}(\epsilon_0-\epsilon)+O(\frac{1}{\tau})\right)\left(\mbox{\boldmath $a$}\cdot\frac{x-p}{\vert x-p\vert}\right)^2.
\end{array}
$$
Since we have 
$$\displaystyle
\left(\mbox{\boldmath $a$}\cdot\frac{x-p}{\vert x-p\vert}\right)^2
=1-\left\vert\mbox{\boldmath $a$}\times\frac{x-p}{\vert x-p\vert}\right\vert^2,
$$
we obtain
$$\begin{array}{l}
\,\,\,\,\,\,
\displaystyle
\frac{\tilde{\epsilon_0}}{\tilde{\epsilon}}(\tilde{\epsilon_0}-\tilde{\epsilon})
\vert\mbox{\boldmath $V$}_e^0\vert^2
+
(\mu-\mu_0)
\vert\mbox{\boldmath $V$}_m^0\vert^2\\
\\
\displaystyle
=\left\{\frac{\epsilon_0}{\epsilon}(\epsilon_0-\epsilon)+
(\mu-\mu_0)\cdot\frac{\epsilon_0}{\mu_0}
+O(\frac{1}{\tau})\right\}K(\tau)^2(\tilde{f}(\tau))^2v^2
\left\vert\mbox{\boldmath $a$}\times\frac{x-p}{\vert x-p\vert}\right\vert^2
\\
\\
\displaystyle
\,\,\
+
\frac{4}{\mu_0\epsilon_0\vert x-p\vert^2}
K(\tau)^2(\tilde{f}(\tau))^2\frac{v^2}{\tau^2}\left(\frac{\epsilon_0}{\epsilon}(\epsilon_0-\epsilon)+O(\frac{1}{\tau})\right).
\end{array}
\tag {3.20}
$$

First consider the case when (B.I) is satisfied.

By Lemmas A.3 in \cite{Ithermo}, we have
$$
\displaystyle
\tau^{3}\int_Dv^2\left\vert\mbox{\boldmath $a$}\times\frac{x-p}{\vert x-p\vert}\right\vert^2\,dx
\ge C e^{-2\tau\sqrt{\mu_0\epsilon_0}\,\text{dist}\,(\{p\},\partial D)}.
\tag {3.21}
$$
Note this is a consequence of the assumption: there exists a point $q\in\partial D\cap S(p,\partial D)$
such that $\partial D$ is locally given by a graph of $C^2$ function on the tangent plane of $\partial D$ at $q$.

We need also a sharp bound of $\Vert v\Vert_{L^2(D)}$.  
Let $q_j$, $j=1,\cdots,m$ be all the points in the first reflector from the point $p$.
By the curvature assumption in statement (ii), for each $j=1,\cdots, m$ one can find an open ball $B'_j$
centered at $q_j$ and a positive number $\lambda_j>0$ such that the set $D\cap B'_j$ is {\it contained} in 
the open ball $B''_j=\{x\in\Bbb R^3\,\vert\,\vert x-(q_j-\lambda_j\mbox{\boldmath $\nu$}_{q_j})\vert<\lambda_j\,\}$. 
And one may choose the balls $B_j'$, $j=1,\cdots,m$ in such a way that $\overline{B'_j}\cap\overline{B'_l}=\emptyset$
if $j\not=j'$ and $\vert x-p\vert\ge \text{dist}\,(\{p\},\partial D)+c$ for all $x\in D\setminus\cup_{j=1}^mB'_j$ with
a positive constant $c$.  Note that $\text{dist}\,(\{p\},\partial D)=\text{dist}\,(\{p\},\partial B_j'')=\vert p-q_j\vert$.

Then one has
$$\begin{array}{ll}
\displaystyle
\int_Dv^2\,dx
&
\displaystyle
=\sum_{j=1}^m\int_{D\cap B'_j}v^2\,dx+\int_{D\setminus\cup_{j=1}^mB'_j}v^2\,dx
\\
\\
\displaystyle
&
\displaystyle
\le \sum_{j=1}^m\int_{B''_j}\frac{e^{-2\tau\,\sqrt{\epsilon_o\mu_0}\,\vert x-p\vert}}{\vert x-p\vert^2}\,dx
+O(e^{-2\tau\,\sqrt{\epsilon_0\,\mu_0}\,\text{dist}\,(\{p\},\partial D)-c'\tau})
\\
\\
\displaystyle
&
\displaystyle
\le C\sum_{j=1}^m\int_{B''_j}\frac{e^{-\sqrt{\epsilon_o\mu_0}\,(2\tau)\,\vert x-p\vert}}{\vert x-p\vert}\,dx
+O(e^{-2\tau\,\sqrt{\epsilon_0\,\mu_0}\,\text{dist}\,(\{p\},\partial D)-c'\tau})
\\
\\
\displaystyle
&
\displaystyle
=C\sum_{j=1}^m
\frac{4\pi\varphi(2\sqrt{\epsilon_0\mu_0}\,\tau\lambda_j)}{(2\sqrt{\epsilon_0\mu_0}\tau)^3}
\,\frac{e^{-2\tau\,\sqrt{\epsilon_0\mu_0}\,\vert p-(q_j-\lambda_j\mbox{\boldmath $\nu$}_{q_j})\vert}}
{\vert p-(q_j-\lambda_j\mbox{\boldmath $\nu$}_{q_j})\vert}
+O(e^{-2\tau\,\sqrt{\epsilon_0\,\mu_0}\,\text{dist}\,(\{p\},\partial D)-c'\tau})
\\
\\
\displaystyle
&
\displaystyle
=O(\tau^{-2})e^{-2\tau\,\sqrt{\epsilon_0\mu_0}\,\text{dist}\,(\{p\},\partial D)},
\end{array}
\tag {3.22}
$$
where $c'=2\,\sqrt{\epsilon_0\mu_0}\,c$.
Note that, at the last step we have applied the mean value theorem \cite{CH} to the integral over $B''_j$
and the equation $\vert p-(q_j-\lambda_j\mbox{\boldmath $\nu$}_{q_j})\vert=\vert p-q_j\vert+\lambda_j$.

By Lemma A.4 in \cite{Ithermo}, we have, as $\tau\rightarrow\infty$
$$\displaystyle
\tau^2\int_Dv^2\,dx\ge C' e^{-2\tau\,\sqrt{\epsilon_0\,\mu_0}\,\text{dist}\,(\{p\},\partial D)}
$$
provided $\partial D$ is locally given by a graph of a $C^2$ function on the tangent plane of $\partial D$ at a point $q\in\partial D\cap S(p,\partial D)$.
So the upper bound (3.22) is sharp.

\noindent
Now applying estimates (3.21) and (3.22) to the integral of right-hand side on (3.20) over $D$ together with (3.15), we obtain (3.16).

Before describing the case when (B.II) is satisfied, it has better to consider the case when (B.III) is satisfied.

From the proof of Lemmas A.3 in \cite{Ithermo}, we have
$$
\displaystyle
\tau^{2}\int_Dv^2\left\vert\mbox{\boldmath $a$}\times\frac{x-p}{\vert x-p\vert}\right\vert^2\,dx
\ge C e^{-2\tau\sqrt{\mu_0\epsilon_0}\,\text{dist}\,(\{p\},\partial D)}.
\tag {3.23}
$$
Note this is a consequence of the assumption: there exists a point $q\in\partial D\cap S(p,\partial D)$
such that $\partial D$ is locally given by a graph of $C^2$ function on the tangent plane of $\partial D$ at $q$
and $\mbox{\boldmath $a$}\times\mbox{\boldmath $\nu$}_q\not=\mbox{\boldmath $0$}$.

On the other hand, we have
$$
\displaystyle
\int_Dv^2\,dx=O(\tau^{-1})e^{-2\tau\,\sqrt{\epsilon_0\mu_0}\,\text{dist}\,(\{p\},\partial D)}.
\tag {3.24}
$$
This is because of:

(a)  the trivial estimate
$$\displaystyle
\int_Dv^2\,dx\le C^{-1}
\int_D\,\frac{e^{-s\vert x-p\vert}}{\vert x-p\vert}\,dx,
\tag {3.25}
$$
where $s=2\sqrt{\epsilon_0\mu_0}\tau$ and $C=\text{dist}\,(\{p\},\partial D)$;

(b) the surface integral expression
$$
\displaystyle
\int_D\,\frac{e^{-s\vert x-p\vert}}{\vert x-p\vert}\,dx
=-s^{-1}\int_{\partial D}
\,e^{-s\vert x-p\vert}\left(1+\frac{1}{s\vert x-p\vert}\,\right)\,\frac{x-p}{\vert x-p\vert^2}\cdot\mbox{\boldmath $\nu_x$}\,dS_x,
\tag {3.26}
$$
which is the consequence the equation
$$\displaystyle
\Delta\left(\frac{e^{-s\vert x-p\vert}}{\vert x-p\vert}\,\right)=s^2\frac{e^{-s\vert x-p\vert}}{\vert x-p\vert}
$$
and integration by parts.

Needless to say, after having (3.23) and (3.24) we obtain the desired conclusion.

Finally consider the case when (B.II) is satisfied.  
In this case one can apply the Laplace method \cite{BH} to the surface integral of the right-hand side on (3.26).
See also \cite{IMax}.
The result is
$$\begin{array}{l}
\displaystyle
\,\,\,\,\,\,
\int_{\partial D}
\,e^{-s\vert x-p\vert}\left(1+\frac{1}{s\vert x-p\vert}\,\right)\,\frac{x-p}{\vert x-p\vert^2}\cdot\mbox{\boldmath $\nu_x$}\,dS_x
\\
\\
\displaystyle
\sim
-s^{-1}
\frac{\pi e^{-s\,\text{dist}\,(\{p\},\partial D)}}{\text{dist}(\{p\},\partial D)}
\sum_{q\in\partial D\cap S(p,\partial D)}
\frac{1}{\displaystyle
\sqrt{\text{det}\,(\mbox{\boldmath $S$}_q(S(p,\partial D))-\mbox{\boldmath $S$}_q(\partial D))}}.
\end{array}
$$
Thus (3.25) and (3.26) yields the estimate
$$\displaystyle
\int_Dv^2\,dx=O(\tau^{-2})e^{-2\tau\,\sqrt{\epsilon_0\mu_0}\,\text{dist}\,(\{p\},\partial D)}.
$$
Now a combination of this and (3.21) yields the desired conclusion.

To obtain (3.18) write
$$\begin{array}{l}
\,\,\,\,\,\,
\displaystyle
(\tilde{\epsilon_0}-\tilde{\epsilon})
\vert\mbox{\boldmath $V$}_e^0\vert^2
+
\frac{\mu_0}{\mu}(\mu-\mu_0)
\vert\mbox{\boldmath $V$}_m^0\vert^2\\
\\
\displaystyle
=\left\{
(\tilde{\epsilon_0}-\tilde{\epsilon})+
\frac{\mu_0}{\mu}(\mu-\mu_0)
\left(\sqrt{\frac{\tilde{\epsilon}_0}{\mu_0}}+\frac{1}{\tau\mu_0\vert x-p\vert}\right)^2
\right\}K(\tau)^2(\tilde{f}(\tau))^2v^2A^2
\left\vert\mbox{\boldmath $a$}\times\frac{x-p}{\vert x-p\vert}\right\vert^2
\\
\\
\displaystyle
\,\,\,
+(\tilde{\epsilon_0}-\tilde{\epsilon})
K(\tau)^2(\tilde{f}(\tau))^2v^2(B-A)^2\left(\mbox{\boldmath $a$}\cdot\frac{x-p}{\vert x-p\vert}\right)^2
\\
\\
\displaystyle
=\left\{
(\epsilon_0-\epsilon)+
\frac{\epsilon_0}{\mu}(\mu-\mu_0)
+O(\frac{1}{\tau})
\right\}
K(\tau)^2(\tilde{f}(\tau))^2v^2A^2
\left\vert\mbox{\boldmath $a$}\times\frac{x-p}{\vert x-p\vert}\right\vert^2
\\
\\
\displaystyle
\,\,\,
+
\frac{4}{\tau^2\mu_0\epsilon_0\vert x-p\vert^2}
K(\tau)^2(\tilde{f}(\tau))^2v^2\left((\epsilon_0-\epsilon)+O(\frac{1}{\tau})\right)\left(\mbox{\boldmath $a$}\cdot\frac{x-p}{\vert x-p\vert}\right)^2
\\
\\
\displaystyle
=\left\{
(\epsilon_0-\epsilon)+
\frac{\epsilon_0}{\mu}(\mu-\mu_0)
+O(\frac{1}{\tau})
\right\}
K(\tau)^2(\tilde{f}(\tau))^2v^2A^2
\left\vert\mbox{\boldmath $a$}\times\frac{x-p}{\vert x-p\vert}\right\vert^2
\\
\\
\displaystyle
\,\,\,
+
\frac{4}{\mu_0\epsilon_0\vert x-p\vert^2}
K(\tau)^2(\tilde{f}(\tau))^2\frac{v^2}{\tau^2}\left((\epsilon_0-\epsilon)+O(\frac{1}{\tau})\right).
\end{array}
$$
Then the remaining parts are similarly done provided (3.17) is valid.

\noindent
$\Box$

{\it\noindent Proof of Lemma 3.3.}
Write
$$
\displaystyle
\frac{\epsilon_0}{\epsilon}(\epsilon_0-\epsilon)+
(\mu-\mu_0)\cdot\frac{\epsilon_0}{\mu_0}
=
\epsilon_0
\left\{
\frac{1}{\epsilon}(\epsilon_0-\epsilon)+
(\mu-\mu_0)\cdot\frac{1}{\mu_0}
\right\}.
$$
We have
$$
\displaystyle
\frac{1}{\epsilon}(\epsilon_0-\epsilon)+
(\mu-\mu_0)\cdot\frac{1}{\mu_0}
=\frac{1}{\epsilon_r}-1+\mu_r-1.
$$
Now the equivalence of (A.I) and (3.15) is clear.

The equivalence of (A.II) and (3.17) is clear since we have
$$\begin{array}{ll}
\displaystyle
(\epsilon_0-\epsilon)+
\frac{\epsilon_0}{\mu}(\mu-\mu_0)
&
\displaystyle
=\epsilon_0
\left\{
\frac{1}{\epsilon_0}(\epsilon_0-\epsilon)-
\frac{1}{\mu}(\mu_0-\mu)
\right\}
\\
\\
\displaystyle
&
\displaystyle
=\epsilon_0
\left(1-\epsilon_r+1-\frac{1}{\mu_r}
\right).
\end{array}
$$

\noindent
$\Box$

\section{Further problems}

Here we point out some of problems to be solved.

$\bullet$  The first problem is to find extraction formulae of the jump of $\mu_r$ and $\epsilon_r$ on $\partial D$.
This is an extraction problem of a {\it quantitative} property of unknown obstacles.
For an impenetrable obstacle with the Leontovich boundary condition, we have already given an extraction formula of
the values of the coefficient in the boundary condition by using the time domain enclosure method, see \cite{IMax3}.
See also \cite{IWd} for a scalar wave case with a dissipative boundary condition on the surface of an obstacle.

$\bullet$  The second one is 
to consider the case when an unknown penetrable obstacle is embedded in a two layered homogeneous background medium.
It will be possible to extend the results in \cite{IK, IKK} which considered a scalar wave equation case, to the Maxwell
system one.

$\bullet$  The third one is the case when the background medium is general inhomogeneous case.
For a scalar wave equation case we have a result in \cite{IWa}.  The proof fully makes use of the
advantage that the governing equation is single.  So it would be interesting to extend the result to the Maxwell system.

$\bullet$  Theorem 1.1 can be extended to the case when the obstacle is inhomogeneous anisotropic and
embedded in an inhomogeneous anisotropic medium.  Then it would be interesting to consider problems 
similar to those mentioned above.

$\bullet$  It would be interesting to consider also find an obstacle embedded in a dispersive metamaterial
in the framework of the direct problem studied in \cite{NV}.  In particular, is there any theorem like Theorem 1.1?

$$\quad$$

\centerline{{\bf Acknowledgment}}

The author was partially supported by Grant-in-Aid for
Scientific Research (C)(No. 17K05331) and (B)(No. 18H01126) of Japan  Society for
the Promotion of Science.

$$\quad$$

\vskip1cm
\noindent
e-mail address

ikehata@hiroshima-u.ac.jp

\end{document}